\documentclass[reqno]{amsart}
\usepackage{amssymb,amsthm,amsfonts,amstext}
\usepackage{amsmath}
\usepackage{graphicx}
\usepackage[american]{babel}
\usepackage[ansinew]{inputenc}
\makeatletter \@addtoreset{equation}{section} \makeatother

\renewcommand\thetable{\thesection.\@arabic\c@table}
\newtheorem{theorem}{Theorem}[section]
\newtheorem{lemma}[theorem]{Lemma}
\newtheorem{proposition}[theorem]{Proposition}
\newtheorem{corollary}[theorem]{Corollary}

\newtheorem{definition}{Definition}
\newtheorem{remark}{Remark}[section]
\newcommand{\dem}{\par\medbreak\noindent{\textbf{Proof of Theorem 2.2.}}}
\newcommand{\demo}{\par\medbreak\noindent{\textbf{Proof of Theorem 2.3.}}}
\newcommand{\demon}{\par\medbreak\noindent{\textbf{Proof of Corollary 2.4.}}}
\newcommand{\demons}{\par\medbreak\noindent{\textbf{Proof of Corollary 2.7}}}
\newcommand{\cqd}{\hfill$\sqcup\!\!\!\!\sqcap\bigskip$}

\title[C.L.T. FOR A TAGGED PARTICLE IN ASYMMETRIC SIMPLE EXCLUSION]{CENTRAL LIMIT THEOREM FOR A TAGGED PARTICLE IN ASYMMETRIC SIMPLE EXCLUSION}
\date{\today}
\begin{document}
\author{Patrícia Gonçalves}
\begin{abstract}
We prove a Functional Central Limit Theorem for the position of a Tagged Particle in the one-dimensional Asymmetric Simple Exclusion Process in
the hyperbolic scaling, starting from a Bernoulli product measure conditioned to have a particle at the origin. We also prove that the position
of the Tagged Particle at time $t$ depends on the initial configuration, by the number of empty sites in the interval $[0,(p-q)\alpha t]$
divided by $\alpha$ in the hyperbolic and in a longer time scale, namely $N^{4/3}$.
\end{abstract}
\renewcommand{\subjclassname}{\textup{2000} Mathematics Subject Classification}
\begin{thanks} {The author wants to express her gratitude to
    F.C.T. (Portugal) for supporting her Phd with the grant /SFRH/ BD/ 11406/ 2002.}
\end{thanks}
\address{IMPA, Estrada Dona Castorina 110, CEP 22460, Rio de Janeiro, Brasil}
\email{patg@impa.br}
\maketitle
\section{Introduction}
The Exclusion process on $\mathbb{Z}^{d}$ has been extensively studied. In this process, particles evolve on $\mathbb{Z}^{d}$ according to
interacting random walks with an exclusion rule which prevents more than one particle per site. The dynamics can be informally described as
follows. Fix a probability $p(\cdot)$ on $\mathbb{Z}^{d}$. Each particle, independently from the others, waits a mean one exponential time, at
the end of which being at $x$ it jumps to $x+y$ at rate $p(y)$. If the site is occupied the jump is suppressed to respect the exclusion rule. In
both cases, the particle waits a new exponential time.

The space state of the process is $\{0,1\}^{\mathbb{Z}^{d}}$ and we denote the configurations by the Greek letter $\eta$, so that $\eta(x)=0$ if
the site $x$ is vacant and $\eta(x)=1$ otherwise. The case in which $p(y)=0$ $\forall |y|>1$ is referred as the Simple Exclusion process and in
the Asymmetric Simple Exclusion process (ASEP) the probability $p$ is such that $p(1)=p$, $p(-1)=1-p$ with $p\neq{1/2}$ while in the Symmetric
Simple Exclusion process (SSEP) $p=1/2$.

For $0\leq{\alpha}\leq{1}$, denote by $\nu_{\alpha}$ the Bernoulli product measure on $\{0,1\}^{\mathbb{Z}^{d}}$ with density $\alpha$. It is
known that $\nu_{\alpha}$ is an invariant measure for the Exclusion process and that all invariant and translation invariant measures are convex
combinations of $\nu_{\alpha}$ if $p(.)$ is such that $p_{t}(x,y)+p_{t}(y,x)>0$, $\forall{x,y\in{\mathbb{Z}^{d}}}$ and $\sum_{x}p(x,y)=1$,
$\forall{y\in{\mathbb{Z}^{d}}}$, see \cite{L.}.

Assume that the origin is occupied at time $0$. Tag this particle and denote by $X_{t}$ its position at time $t$. Applying an invariance
principle due to Newman and Wright \cite{N.}, Kipnis in \cite{K.} proved a C.L.T. for the position of the Tagged Particle in the one-dimensional
ASEP, provided the initial configuration is distributed according to $\nu_{\alpha}^{*}$, the Bernoulli product measure conditioned to have a
particle at the origin. Transforming the exclusion process into a series of queues, an asymmetric Zero-Range process with constant rate, the
position of the Tagged Particle becomes the current through the bond $[-1,0]$. Kipnis \cite{K.}, was able to apply Newman and Wright results to
the Zero-Range process and derive the L.L.N. and C.L.T. for the position of the Tagged Particle.

Few years later, Ferrari and Fontes \cite{F.F.2} proved that the position at time t of the Tagged Particle, $X_{t}$, can be approximated by a
Poisson Process. More precisely, they proved that for all $t\geq{0}$, if the initial distribution is $\nu_{\alpha}^{*}$ and $p>q$,
$X_{t}=N_{t}-B_{t}+B_{0}$, where $N_{t}$ is a Poisson Process with rate $(p-q)(1-\alpha)$ and $B_{t}$ is a stationary process with bounded
exponential moments. As a corollary they obtained the weak convergence of
\begin{equation*}
\frac{X_{t{\epsilon}^{-1}}-(p-q)(1-\alpha)t{\epsilon}^{-1}}{\sqrt{(p-q)(1-\alpha)t{\epsilon}^{-1}}}
\end{equation*}
to a Brownian motion. The argument is divided in two steps. The convergence of the finite-dimensional distributions \cite{F.} is consequence of
the fact that in the scale $t^{\frac{1}{2}}$, the position $X_{t}$ can be read from the initial configuration: $X_{t}$ is given by the initial
number of empty sites in the interval $[0,(p-q)\alpha t]$ divided by $\alpha$. Tightness follows from the sharp approximation of
$X_{t{\epsilon}^{-1}}$ by the Poisson process and the weak convergence of the Poisson process to Brownian motion. Using the approximation of
$X_{t}$ by a Poisson process and Kipnis results for the Tagged Particle, the same authors prove equilibrium density fluctuations for the ASEP in
\cite{F.F.1}. The density fluctuations for the Totally Asymmetric Simple Exclusion process (the case $p=1$) have also been obtained by
Rezakhanlou in \cite{Reza.} in a more general setting than for the process starting from an equilibrium state.

Recently, Jara and Landim in \cite{J.L.} showed that the asymptotic behavior of the Tagged Particle in the one-dimensional nearest neighbor
exclusion process, can be recovered from a joint asymptotic behavior of the empirical measure and the current through a bond. From this
observation they proved a non-equilibrium C.L.T. for the position of the Tagged Particle in the SSEP, under diffusive scaling.

In this paper, besides using this general method to reprove Ferrari and Fontes result on the convergence of the rescaled position of the Tagged
Particle to a Brownian motion in the hyperbolic time scale, we extended their result by showing that in a longer time scale the position of the
Tagged Particle still depends on the initial configuration.

The advantage of our approach is that it relates the C.L.T. for the position of the Tagged Particle to the C.L.T. for the empirical measure, a
problem which is relatively well understood, see \cite{K.L.}. In particular, we can expect to apply this approach for a one-dimensional system
in contact with reservoirs.

It was shown by Rezakhanlou in \cite{R.}, that in the ASEP the macroscopic particle density profile in the hyperbolic scaling, evolves according
to the inviscid Burgers equation, namely: $\partial_{t}\rho(t,u)+(p-q)\nabla(\rho(t,u)(1-\rho(t,u)))=0$. To establish the C.L.T. for the
empirical measure we need to consider the density fluctuation field as defined in (\ref{eq:densfieldinz}) below. We show that, in this time
scale, the time evolution of the limit density fluctuation field is deterministic, in the sense that at any given time $t$, the density field is
a translation of the initial one. As mentioned above, this result was previously obtained in \cite{F.F.1}. In order to observe fluctuation from
the dynamics one has to change to the diffusive scaling (see section six).

The translation or velocity of the system is given by $v=(p-q)(1-2\alpha)$ and for $\alpha=1/2$, the field does not evolve in time, and one is
forced to go beyond the hydrodynamic scaling. We can consider the density fluctuation field in the longer time scale as defined in
(\ref{eq:densfieldlongscale}), where we subtract the velocity of the system and any value of $\alpha$ can be considered in this setting.

It is conjectured that until the time scale $N^{3/2}$ the density fluctuation field does not evolve in time, see chap.5 of \cite{HS.} and
references therein. The result we obtain is a contribution in this direction, since we can accomplish the result just up to the time scale
$N^{4/3}$. The main difficulty in proving the C.L.T. for the empirical measure is the Boltzmann-Gibbs Principle, which we are able to prove for
this time scale using a multi-scale argument.

As a consequence of this translation behavior, we show the dependence on the initial configuration of the current through a bond and the
position of the Tagged Particle in the longer time scale.

This work is organized as follows. In the second section we introduce some notation and we state the results. The sketch of the proof of the
C.L.T. for the empirical measure associated to the ASEP in the hyperbolic scaling is exposed in the third section. In section four, we use the
same strategy as in \cite{J.L.} to obtain L.L.N. and the convergence of finite-dimensional distributions of the position of the Tagged Particle
to those of Brownian motion. Tightness is proved by means of the Zero-Range representation as Kipnis in \cite{K.}. In this time scale we show
that the current through a fixed bond and the position of the Tagged Particle at time $tN$, can be read from the initial configuration, in
section five.

In the following sections we study the same problem up to the time scale $N^{1+\gamma}$ with ${\gamma<1/3}$.  We start by showing the C.L.T. for
the empirical measure associated to this process, in section six. Since a Boltzmann-Gibbs Principle is needed, its proof is the content of the
seventh section and in the subsequent section we treat the problem of tightness. In the last section we prove the dependence on the initial
configuration for the current through a bond that depends on time and the position of the Tagged Particle, in this longer time scale.
\section{Statement of Results}
The one-dimensional asymmetric simple exclusion process is the Markov Process $\eta_{t}\in{\{0,1\}^{\mathbb{Z}}}$ with generator given on local
functions by
\begin{equation}
Lf(\eta)=\sum_{x\in{\mathbb{Z}}}\sum_{y=x\pm{1}}c(x,y,\eta)[f(\eta^{x,y})-f(\eta)], \label{eq:generator}
\end{equation}
where $c(x,y,\eta)=p(x,y)\eta(x)(1-\eta(y))$, $p(x,x+1)=p$, $p(x,x-1)=q=1-p$ and
\[ \eta^{x,y}(z)=\left\{
\begin{array}{rl}
\eta(z), & \mbox{if $z\neq{x,y}$}\\
\eta(y), & \mbox{if $z=x$}\\
\eta(x), & \mbox{if $z=y$}
\end{array}.
\right.
\]
\\
Its description is the following. At most one particle is allowed at each site. If there is a particle at site $x$, it jumps at rate $p$ to site
$x+1$ if there is no particle at that site. If the site $x-1$ is empty, the particle at $x$ jumps to $x-1$ at rate $q$. Initially, place the
particles according to a Bernoulli product measure in $\{0,1\}^{\mathbb{Z}}$, of parameter $\alpha\in(0,1)$, denoted by $\nu_{\alpha}$.

For each configuration $\eta$, denote by $\pi^{N}(\eta,du)$ the empirical measure given by
\begin{equation*}
\pi^{N}(\eta,du)=\frac{1}{N}\sum_{x\in\mathbb{Z}}\eta(x)\delta_{\frac{x}{N}}(du)
\end{equation*} and let $\pi_{t}^{N}(\eta,du)=\pi^{N}(\eta_{t},du)$.
First, we state the C.L.T. for the empirical measure, for which we need to introduce some notation.

For each integer $z\geq{0}$, let $H_{z}(x)=(-1)^{z}e^{x^{2}}\frac{d^{z}}{dx}e^{-x^{2}}$ be the Hermite polynomial, and
$h_{z}(x)=\frac{1}{c_{z}}H_{z}(x)e^{-x^{2}}$ the Hermite function, where $c_{z}=z!\sqrt{2\pi}$ . The set $\{h_{z},z\geq{0}\}$ is an orthonormal
basis of $L^{2}(\mathbb{R})$. Consider in $L^{2}(\mathbb{R})$ the operator $K_{0}=x^{2}-\Delta$. A simple computation shows that
$K_{0}h_{z}=\gamma_{z}h_{z}$ where $\gamma_{z}=2z+1$.

For an integer $k\geq{0}$, denote by $\mathcal{H}_{k}$ the Hilbert space induced by $S(\mathbb{R})$ (the space of smooth rapidly decreasing
functions) and the scalar product $<\cdot,\cdot>_{k}$ defined by $<f,g>_{k}=<f,K_{0}^{k}g>$, where $<\cdot,\cdot>$ denotes the inner product of
$L^{2}(\mathbb{R})$ and denote by $\mathcal{H}_{-k}$ the dual of $\mathcal{H}_{k}$, relatively to this inner product.

Fix $\alpha\in{(0,1)}$ and an integer $k$. Denote by $Y_{.}^{N}$ the density fluctuation field, a linear functional acting on functions $H\in
S(\mathbb{R})$ as
\begin{equation}
Y_{t}^{N}(H)=\frac{1}{\sqrt{N}}\sum_{x\in{\mathbb{Z}}}H\Big(\frac{x}{N}\Big)(\eta_{tN}(x)-\alpha). \label{eq:densfieldinz}
\end{equation}

Denote by $D(\mathbb{R}^{+},\mathcal{H}_{-k})$ (resp. $C(\mathbb{R}^{+},\mathcal{H}_{-k})$) the space of $H_{-k}$-valued functions, right
continuous with left limits (resp. continuous), endowed with the uniform weak topology, by $Q_{N}$ the probability measure on
$D(\mathbb{R}^{+},\mathcal{H}_{-k})$ induced by the density fluctuation field $Y^{N}_{.}$ and $\nu_{\alpha}$. Consider
$\mathbb{P}^{N}_{\nu_{\alpha}}=\mathbb{P}_{\nu_{\alpha}}$ the probability measure on $D(\mathbb{R}^{+},\{0,1\}^{\mathbb{Z}})$ induced by
$\nu_{\alpha}$ and the Markov process $\eta_{t}$ speeded up by $N$ and denote by $\mathbb{E}_{\nu_{\alpha}}$ expectation with respect to
$\mathbb{P}_{\nu_{\alpha}}$.

\begin{theorem} \label{th:flu1}
Fix an integer $k>2$. Denote by $Q$ be the probability measure on $C(\mathbb{R}^{+},\mathcal{H}_{-k})$ corresponding to a stationary Gaussian
process with mean $0$ and covariance given by
\begin{equation}
E_{Q}[Y_{t}(H)Y_{s}(G)]=\chi(\alpha)\int_{\mathbb{R}}H(u+v(t-s))G(u)du \label{eq:covar}
\end{equation}
for every $0\leq{s}\leq{t}$ and $H$, $G$ in $\mathcal{H}_{k}$. Here $\chi(\alpha)= \textbf{Var}(\nu_{\alpha},\eta(0))$ and
$v=(p-q)\chi'(\alpha)$. Then, the sequence $(Q_{N})_{N\geq{1}}$ converges weakly to the probability measure $Q$.
\end{theorem}
We remark that last theorem holds for the ASEP evolving in any $\mathbb{Z}^{d}$, with the appropriate changes. In this case, the limit density
fluctuation field at time $t$ is a translation of the initial density field, since for every $H\in{S(\mathbb{R})}: Y_{t}(H)=Y_{0}(T_{t}H)$,
where $T_{t}H(u)=H(u+vt)$.

Having established the equilibrium density fluctuations, we can obtain the L.L.N. and the C.L.T. for the current over a bond, as in \cite{J.L.}.
Denote by $\nu_{\alpha}^{*}$ the measure $\nu_{\alpha}$ conditioned to have a particle at the origin. By coupling the ASEP starting from
$\nu_{\alpha}$ with the ASEP starting from $\nu_{\alpha}^{*}$, in such a way that both processes differ at most in one site at any given time,
the L.L.N. and the C.L.T. for the empirical measure and for the current starting from $\nu_{\alpha}^{*}$, follows from the L.L.N. and the C.L.T.
for the empirical measure and for the current starting from $\nu_{\alpha}$.

Assume now that the initial measure is $\nu_{\alpha}^{*}$, let $\mathbb{P}^{N}_{\nu_{\alpha}^{*}}=\mathbb{P}_{\nu_{\alpha}^{*}}$ be the
probability measure on $D(\mathbb{R}^{+},\{0,1\}^{\mathbb{Z}})$ induced by $\nu_{\alpha}^{*}$ and the Markov process $\eta_{t}$ speeded up by
$N$ and denote by $\mathbb{E}_{\nu_{\alpha}^{*}}$ expectation with respect to $\mathbb{P}_{\nu_{\alpha}^{*}}$.

Denote by $X_{tN}$ the position at time $tN\geq{0}$ of the tagged particle initially at the origin. We reprove the L.L.N. for the position of
the Tagged Particle, which was previously obtained by Saada in \cite{E.S.}:

\begin{theorem} \label{th:llntagged}
Fix $t\geq{0}$. Then,
\begin{equation*}
\frac{X_{tN}}{N}\xrightarrow[N\rightarrow{+\infty}]\,v_{t}=(p-q)(1-\alpha)t
\end{equation*}
in $\mathbb{P}_{\nu_{\alpha}^{*}}$-probability.
\end{theorem}
and the convergence to the Brownian motion, which was already obtained by Ferrari and Fontes in \cite{F.F.2}:
\begin{theorem}\label{th:clttagged}
Under $\mathbb{P}_{\nu_{\alpha}^{*}}$,
\begin{equation*}
\frac{X_{tN}-v_{t}N}{\sqrt{N|p-q|(1-\alpha)}}\xrightarrow[N\rightarrow{+\infty}]\,B_{t}
\end{equation*}
weakly, where $B_{t}$ denotes the standard Brownian motion.
\end{theorem}

Another interesting property is the dependence on the initial configuration for the position of the Tagged Particle, which was previously
obtained by Ferrari in \cite{F.}. Suppose $p>q$.
\begin{corollary} \label{tagged dependence}
Fix $t\geq{0}$. Then for every $\epsilon>0$,
\begin{equation*}
\lim_{N\rightarrow{+\infty}}\mathbb{E}_{\nu_{\alpha}^{*}}\Big[\frac{{X_{tN}}}{\sqrt{N}}-\frac{\sum_{x=0}^{(p-q)\alpha tN}(1-\eta_{0}(x))}{\alpha
\sqrt{N}}\Big]^{2-\epsilon}=0.
\end{equation*}
\end{corollary}

In the hyperbolic scaling, we have seen above that for the case $\alpha=1/2$ the limit density fluctuation field at time $t$ is the same as the
initial one. This forced us to consider a longer time scale in order to observe other fluctuations than the shifted version of the initial ones.

Henceforth, consider the ASEP evolving in the time scale $N^{1+\gamma}$, with $\gamma>0$. In the sequel, we point out the restrictions needed in
$\gamma$ in order to obtain the results.

Let $\alpha\in(0,1)$ and redefine the density fluctuation field on $H\in{S(\mathbb{R})}$ by:
\begin{equation}
Y_{t}^{N,\gamma}(H)=\frac{1}{\sqrt{N}}\sum_{x\in{\mathbb{Z}}}H\Big(\frac{x-vtN^{1+\gamma}}{N}\Big)(\eta_{t{N}^{1+\gamma}}(x)-\alpha).
\label{eq:densfieldlongscale}
\end{equation}

We remark here, than one can define in the hyperbolic scaling of time the density fluctuation field as above. But in that case the current could
not be defined through a fixed bond, instead it would have to be defined through a bond that depends on time (see section 9). As we we want to
show the C.L.T. for the position of a Tagged Particle using the relation between the density of particles and the current through a fixed bond
(\ref{eq:taggedparticle}), we have the need to defined the density fluctuation field as in (\ref{eq:densfieldinz}).

As above, let $Q^\gamma_{N}$ be the probability measure on $D(\mathbb{R}^{+},\mathcal{H}_{-k})$ induced by the density fluctuation field
$Y^{N,\gamma}_{.}$ and $\nu_{\alpha}$, let $\mathbb{P}^{N,\gamma}_{\nu_{\alpha}}=\mathbb{P}^{\gamma}_{\nu_{\alpha}}$ be the probability measure
on $D(\mathbb{R}^{+},\{0,1\}^{\mathbb{Z}})$ induced by $\nu_{\alpha}$ and the Markov process $\eta_{t}$ speeded up by $N^{1+\gamma}$ and denote
by $\mathbb{E}_{\nu_{\alpha}}^{\gamma}$ expectation with respect to $\mathbb{P}_{\nu_{\alpha}}^{\gamma}$. Now, we state Theorem \ref{th:flu1} in
this longer scaling:

\begin{theorem}
Fix an integer $k>1$ and $\gamma<1/3$. Let $Q$ be the probability measure on $C(\mathbb{R}^{+},\mathcal{H}_{-k})$ corresponding to a stationary
Gaussian process with mean $0$ and covariance given by
\begin{equation}\label{eq:covariancelonger}
E_{Q}[Y_{t}(H)Y_{s}(G)]=\chi(\alpha)\int_{\mathbb{R}}H(u)G(u)du
\end{equation}
for every $s,t\geq{0}$ and $H$, $G$ in $\mathcal{H}_{k}$. Then, the sequence $(Q^{\gamma}_{N})_{N\geq{1}}$ converges weakly to the probability
measure $Q$. \label{th:flu2}
\end{theorem}

As we follow the martingale approach, the main difficulty in proving this theorem is the Boltzmann-Gibbs Principle, which we can prove for
$\gamma<1/3$ and in this case is stated in the following way:

\begin{theorem}{(Boltzmann-Gibbs Principle)} \label{bg}

Fix $\gamma<1/3$. For every $t>0$ and $H\in{S(\mathbb{R})}$,
\begin{equation*}
\lim_{N\rightarrow{\infty}}\mathbb{E}_{\nu_{\alpha}}^{\gamma}\Big[\int_{0}^{t}
\frac{N^{\gamma}}{\sqrt{N}}\sum_{x\in{\mathbb{Z}}}H\Big(\frac{x}{N}\Big)\bar{\eta}_{s}(x)\bar{\eta}_{s}(x+1)ds\Big]^{2}=0.
\end{equation*}
\end{theorem}

In order to keep notation simple, here and after we denote by $\bar{X}$ the centered random variable $X$. Let
$\mathbb{P}^{N,\gamma}_{\nu_{\alpha}^{*}}=\mathbb{P}^{\gamma}_{\nu_{\alpha}^{*}}$ be the probability measure on
$D(\mathbb{R}^{+},\{0,1\}^{\mathbb{Z}})$ induced by $\nu_{\alpha}^{*}$ and the Markov process $\eta_{t}$ speeded up by $N^{1+\gamma}$.

By the results just stated, in this longer time scale the system translates in time at a certain velocity $v$. This allows us to deduce from the
previous results the asymptotic behavior of the position of the Tagged Particle even in the longer time scale:

\begin{corollary} \label{tagged dependence higher}
Fix $t\geq{0}$, suppose that $p>q$ and $\gamma<1/3$. Then,

\begin{equation*}
\frac{{X_{tN^{1+\gamma}}}}{\sqrt{N}}-\frac{\sum_{x=0}^{(p-q)\alpha tN^{1+\gamma}}(1-\eta_{0}(x))}{\alpha
\sqrt{N}}\xrightarrow[N\rightarrow{+\infty}]\,0
\end{equation*}
in $\mathbb{P}_{\nu_{\alpha}^{*}}^{\gamma}$-probability.
\end{corollary}

\section{Density Fluctuations in the Hyperbolic Scaling}
The aim of this section is to prove Theorem \ref{th:flu1}. We just give a sketch of the proof of this result, since we are going to use similar
techniques to the ones used in chap. 11 of \cite{K.L.} when describing the equilibrium fluctuation field of the Symmetric Zero-Range process
under diffusive scaling.

Fix a positive integer $k$ and recall the definition of the density fluctuation field in (\ref{eq:densfieldinz}). The purpose is to show that
$Y_{.}^{N}$ converges to a process $Y_{.}$ whose time-evolution is deterministic.

Denote by $\mathfrak{A}$ the operator $v\nabla$ defined on a domain of $L^{2}(\mathbb{R})$ and by $\{T_{t}, t\geq0\}$ the semigroup associated
to $\mathfrak{A}$. For $t\geq0$, let $\mathcal{F}_{t}$ be the $\sigma$-algebra on $D([0,T],\mathcal{H}_{-k})$ generated by $Y_{s}(H)$ for
$s\leq{t}$ and $H$ in $S(\mathbb{R})$ and set $\mathcal{F}=\sigma(\bigcup_{t\geq{0}}\mathcal{F}_{t})$.

To prove the theorem we need to verify that $(Q_{N})_{N\geq{1}}$ is tight and to characterize the limit field. To check the last assertion, we
consider a collection of martingales associated to the empirical measure. Fix a function $H\in{S(\mathbb{R})}$. Then:
\begin{equation*} \label{eq:martingale M}
M^{N,H}_{t}=Y^{N}_{t}(H)-Y^{N}_{0}(H)-\int^{t}_{0}\frac{1}{\sqrt{N}}\sum_{x\in\mathbb{Z}}\nabla^{N}{H\Big(\frac{x}{N}\Big)}W_{x,x+1}(\eta_{s})ds
\end{equation*}
is a martingale with respect to the filtration $\tilde{\mathcal{F}}_{t}=\sigma(\eta_{s}, s\leq{t})$, whose quadratic variation is given by:
\begin{equation*} \label{eq:quadraticvariation}
\int^{t}_{0}\frac{1}{N^2}\sum_{x\in\mathbb{Z}}\Big(\nabla^{N}{H\Big(\frac{x}{N}\Big)}\Big)^2\Big[c(x,x+1,\eta_{s})+c(x+1,x,\eta_{s})\Big]ds,
\end{equation*}
where $W_{x,x+1}(\eta)$ denotes the instantaneous current between the sites $x$ and $x+1$:
\begin{equation*}
W_{x,x+1}(\eta)=c(x,x+1,\eta)-c(x+1,x,\eta)
\end{equation*}
and
\begin{equation*}
\nabla^{N}H\Big(\frac{x}{N}\Big)=N\Big(H\Big(\frac{x+1}{N}\Big)-H\Big(\frac{x}{N}\Big)\Big).
\end{equation*}

Using the fact that $\sum_{x\in{\mathbb{Z}}}\nabla^{N}{H(\frac{x}{N})}=0$, the integral part of the martingale is equal to:
\begin{equation*}
\int^{t}_{0}\frac{1}{\sqrt{N}}\sum_{x\in\mathbb{Z}}\nabla^{N}{H\Big(\frac{x}{N}\Big)} \Big[\bar{W}_{x,x+1}(\eta_{s})\Big]ds.
\end{equation*}
As we need to write the expression inside last integral in terms of the fluctuation field $Y_{s}^{N}$, we are able to replace the function
$\bar{W}_{x,x+1}(\eta_{s})$ by $(p-q)\chi'(\alpha)[\eta_{s}(x)-\alpha]$, with the use of the:

\begin{theorem}{(Boltzmann-Gibbs Principle)} \label{th:bg}

For every local function $g$, for every $H\in{S(\mathbb{R})}$ and every $t>0$,
\begin{equation*}
\lim_{N\rightarrow{\infty}}\mathbb{E}_{\nu_{\alpha}}\Big[\Big(\int_{0}^{t}
\frac{1}{\sqrt{N}}\sum_{x\in{\mathbb{Z}}}H\Big(\frac{x}{N}\Big)\Big\{\tau_{x}g(\eta_{s})
-\tilde{g}(\alpha)-\tilde{g}'(\alpha)[\eta_{s}(x)-\alpha]\Big\}ds\Big)^2\Big]=0,
\end{equation*}
where $\tilde{g}(\alpha)=E_{\nu_{\alpha}}[g(\eta)]$.
\end{theorem}

In spite of considering the ASEP in the hyperbolic scaling, the proof of last result is very close to the one presented for the Zero-Range
process in the diffusive scaling, and for that reason we have omitted it.

Assume now, that $(Q_{N})_{N\geq{1}}$ is tight and let $Q$ be one of its limiting points. By the result just stated and since
$\lim_{N\rightarrow{+\infty}}\mathbb{E}_{\nu_{\alpha}}[(M_{t}^{N,H})^{2}]=0$, under $Q$
\begin{equation} \label{timeevolutionofdfield}
Y_{t}(H)=Y_{0}(H)+\int_{0}^{t}Y_{s}(\mathfrak{A}H)ds.
\end{equation}
 So, $\frac{d}{dt}Y_{t}(H)=Y_{t}(\mathfrak{A}H)$. Take
$r<t$, and note that $\frac{d}{dr}<Y_{r},T_{t-r}H>=0$. As a consequence, $Y_{t}(H)=Y_{0}(T_{t}H)$ where $T_{t}H(u)=H(u+vt)$.

It is easy to show that $Q$ restricted to $\mathcal{F}_{0}$, is a Gaussian field with covariance given by
$E_{Q}(Y_{0}(G)Y_{0}(H))=\chi(\alpha)<G,H>$ and it is immediate that the limit field has covariance given by (\ref{eq:covar}).

To finish the proof, it remains to show that $(Q_{N})_{N\geq{1}}$ is tight whose proof follows closely the same arguments as the ones for the
Zero-Range process in the diffusive scaling. Lastly, we note that once the process evolves on $\mathbb{Z}$ and the hyperbolic scale is
considered, we must take an integer $k>2$ in order have the density fluctuations field well defined in $\mathcal{H}_{-k}$.

\section{Law of Large Numbers and Central Limit Theorem for the Position of the Tagged Particle}

In this section we prove Theorems \ref{th:llntagged} and 2.3 following the same arguments as Jara and Landim in \cite{J.L.}. For that reason we
give an outline of the proofs.

First we state the C.L.T. for the current through a fixed bond. For a site $x$, denote the current through the bond $[x,x+1]$ by
$J^N_{x,x+1}(t)$, as the total number of jumps from the site $x$ to the site $x+1$ minus the total number of jumps from the site $x+1$ to the
site $x$ during the time interval $[0,tN]$. Since
\begin{equation*}
J^N_{-1,0}(t)=\sum_{x\geq{0}}\Big(\eta_{t}(x)-\eta_{0}(x)\Big),
\end{equation*}
the current can be written in terms of the density fluctuation field as
\begin{equation*}
\frac{1}{\sqrt{N}}\Big\{J^N_{-1,0}(t)-\mathbb{E}_{\nu_{\alpha}}[J^N_{-1,0}(t)]\Big\}=Y_{t}^{N}(T_{t}H_{0})-Y_{0}^{N}(H_{0}),
\end{equation*}
where $H_{0}$ is the Heaviside function, $H_{0}(u)=1_{[0,\infty)}(u)$. By approximating $H_{0}$ by a sequence $(G_{n})_{n\geq{1}}$, defined for
each $u\in{\mathbb{R}}$ by $G_{n}(u)=(1-\frac{u}{n})^{+}1_{[0,\infty)}(u)$, we obtain:

\begin{proposition} \label{prop:1}
For every $t\geq{0}$,
\begin{equation*}
\lim_{n\rightarrow{+\infty}}\mathbb{E}_{\nu_{\alpha}}\Big[\frac{\bar{J}^N_{-1,0}(t)}{\sqrt{N}}-(Y_{t}^{N}(T_{t}G_{n})-Y_{0}^{N}(G_{n}))\Big]^{2}=0
\end{equation*}
uniformly in $N$.
\end{proposition}
\begin{proof}
For a site $x$, consider the martingale $M^N_{x,x+1}(t)$ equal to
\begin{equation}
J^N_{x,x+1}(t)-\int_{0}^{t}N W_{x,x+1}(\eta_{s})-vN\eta_{s}(x)ds \label{martingaledecompofcurrent}
\end{equation}
whose quadratic variation is given by
\begin{equation*}
<M^N_{x,x+1}>_{t}=N\int_{0}^{t} \Big\{c(x,x+1,\eta_{s})+c(x+1,x,\eta_{s})\Big\}ds.
\end{equation*}
Since the number of particles is preserved, it holds that:
\begin{equation*}
J^N_{x-1,x}(t)-J^N_{x,x+1}(t)=\eta_{t}(x)-\eta_{0}(x)
\end{equation*}
for all $x\in{\mathbb{Z}}$, $t\geq{0}$, and we have that
\begin{equation*}
Y_{t}^{N}(T_{t}G_{n})-Y_{0}^{N}(G_{n})=\frac{1}{\sqrt{N}}\sum_{x\in{\mathbb{Z}}}G_{n}\Big(\frac{x}{N}\Big)\Big\{\bar{J}^N_{x-1,x}(t)-\bar{J}^N_{x,x+1}(t)\Big\}.
\end{equation*}
Making a summation by parts and using the explicit knowledge of $G_{n}$, last expression can be written as
\begin{equation*}
\frac{\bar{J}^N_{-1,0}(t)}{\sqrt{N}}-\Big[Y_{t}^{N}(T_{t}G_{n})-Y_{0}^{N}(G_{n})\Big]=\frac{1}{\sqrt{N}}\sum_{x=1}^{Nn}\frac{1}{Nn}\bar{J}^N_{x-1,x}(t).
\end{equation*}
Representing the current ${J}^N_{x-1,x}(t)$ in terms of the martingales $M^N_{x-1,x}(t)$, the right hand side of the last expression becomes
equal to
\begin{equation} \label{eq:current}
\frac{1}{\sqrt{N}}\sum_{x=1}^{Nn}\frac{1}{Nn}M^N_{x-1,x}(t)
+\frac{1}{\sqrt{N}}\int_{0}^{t}\frac{1}{n}\sum_{x=1}^{Nn}[\bar{W}_{x-1,x}(\eta_{s})-v(\eta_{s}(x-1)-\alpha)]ds.
\end{equation}
The martingale term converges to $0$ in $L^{2}(\mathbb{P}_{\nu_{\alpha}})$ as $n\rightarrow{+\infty}$, since we can estimate their quadratic
variation by $Nt$, use the fact that they are orthogonal to obtain that its $L^{2}(\mathbb{P}_{\nu_{\alpha}})$-norm is bounded above by
$\frac{Ct}{Nn}$.

Making an elementary computation it is easy to show that the $L^{2}(\mathbb{P}_{\nu_{\alpha}})$-norm of the integral term is bounded above by
$\frac{C}{n}$. Taking the limit as $n\rightarrow{\infty}$, the proof is concluded.
\end{proof}
Putting together, last result and the C.L.T. for the empirical measure, it holds:
\begin{theorem} \label{th:cltcurrent}
Fix $x\in{\mathbb{Z}}$ and let
\begin{equation*}
Z_{t}^{N}=\frac{1}{\sqrt{N}}\Big\{J^N_{x,x+1}(t)-\mathbb{E}_{\nu_{\alpha}}[J^N_{x,x+1}(t)]\Big\}.
\end{equation*}
Then, for every $k\geq{1}$ and every $0\leq{t_{1}}<{t_{2}}<..<t_{k}$, $(Z_{t_{1}}^{N},..,Z_{t_{k}}^{N})$ converges in law to a Gaussian vector
$(Z_{t_{1}},..,Z_{t_{k}})$ with mean zero and covariance given by
\begin{equation*}
E_{Q}[Z_{t}Z_{s}]=\chi(\alpha)|v|s
\end{equation*}
provided $s\leq{t}$.
\end{theorem}

Assume now, the initial measure to be $\nu_{\alpha}^{*}$. Let $X_{tN}$ be the position of the Tagged Particle at time $tN\geq{0}$ initially at
the origin. Fix, a positive integer $n$. Since we are considering the one-dimensional setting, particles cannot jump over other particles, and
therefore it holds the following relation:
\begin{equation}
\{X_{tN}\geq{n}\}=\Big\{J^N_{-1,0}(t)\geq{\sum_{x=0}^{n-1}\eta_{t}(x)}\Big\} \label{eq:taggedparticle}
\end{equation}
which allows, together with the previous results, to obtain L.L.N. and the C.L.T. for the position of the Tagged Particle. Now, we give a sketch
of the proof of this results.

\quad\

 \dem \\ In order to show L.L.N. for the Tagged Particle, denote by $\lceil{a}\rceil$ the smallest integer larger or equal to $a$, fix
$u>0$ and take $n=\lceil{uN}\rceil$ in (\ref{eq:taggedparticle}). Using the martingale decomposition of the current
(\ref{martingaledecompofcurrent}) and Theorem \ref{th:cltcurrent} it is easy to show that
\begin{equation*}
\frac{J^N_{-1,0}(t)}{N}\xrightarrow[N\rightarrow{+\infty}]\,(p-q)\chi(\alpha)t
\end{equation*}
in $\mathbb{P}_{\nu_{\alpha}}$-probability. Since $<\pi_{t}^{N},1_{[0,u]}>$ converges in probability to $\alpha u$, we obtain that

\[ \lim_{N\rightarrow{+\infty}}\mathbb{P}_{\nu_{\alpha}^{*}}\Big[\frac{X_{tN}}{N}\geq{u}\Big]=\left\{
\begin{array}{rl}
0, & \mbox{if $(p-q)\chi(\alpha)t<\alpha u$}\\
1, & \mbox{if $(p-q)\chi(\alpha)t\geq{\alpha u}$}
\end{array}.
\right.
\]
\\
For $u<0$ we obtain a similar result, which concludes the proof. \cqd

We proceed by proving the convergence of the Tagged Particle process, properly centered and rescaled, to the standard Brownian motion.

\quad\

\demo \\Let $W_{tN}=\frac{1}{\sqrt{N}}(X_{tN}-v_{t}N)$. The result follows from showing the convergence of finite dimensional distributions of
$W_{tN}$ to those of Brownian motion together with tightness.

Using (\ref{eq:taggedparticle}), Theorems \ref{th:flu1} and \ref{th:cltcurrent} above, it is not hard to show that under
$\mathbb{P}_{\nu_{\alpha}^{*}}$, $\forall{k\geq{1}}$, $\forall{0\leq{t_{1}}<..<t_{k}}$, $(W_{t_{1}N},..,W_{t_{k}N})$ converges in law to a
Gaussian vector $(W_{t_{1}},...,W_{t_{k}})$ with mean zero and covariance given by
\begin{equation*}
E_{Q}\Big[W_{t}W_{s}\Big]=|p-q|(1-\alpha)s
\end{equation*}
for $0\leq{s}\leq{t}$.

To end the proof it remains to show tightness. For that we use a relation between the ASEP and a Zero-Range process, as Kipnis in \cite{K.}. For
the latter, the product measures $\mu_{\alpha}$ with marginals given by $\mu_{\alpha}\{\eta(x)=k\}=\alpha(1-\alpha)^{k}$ are extremal invariant.

This process has space state $\mathcal{X}=\mathbb{N}^{\mathbb{Z}}$ and generator defined on local functions by
\begin{equation*}
\Omega f(\eta)=\sum_{x\in{\mathbb{Z}}}1_{\{\eta(x)\geq{1}\}}[pf(\eta^{x,x-1})+qf(\eta^{x,x+1})-f(\eta)],
\end{equation*}
where $p+q=1$ and
\[\eta^{x,y}(z)=
\begin{cases}
\eta(z), & \mbox{if $z\neq{x,y}$}\\
\eta(x)-1, & \mbox{if $z=x$}\\
\eta(y)+1, & \mbox{if $z=y$}
\end{cases}.
\]
\\
The process can also be reversed with respect to any $\mu_{\alpha}$, and the reversed process is denoted by $\hat{\eta}$, whose generator
$\hat{\Omega}$ is the same as $\Omega$, except that p is replaced by $q$ and vise-versa.

The position of the Tagged Particle in the Zero-Range representation becomes the current through the bond $[-1,0]$:
$X_{t}=-N_{t}^{+}+N_{t}^{-}$, where $N_{t}^{+}$ (resp. $N_{t}^{-}$) is the number of particles that jumped from site $-1$ to site $0$ during the
time interval $[0,t]$ (resp. from site $0$ to $-1$).

As a consequence, the proof ends if we show tightness of the distributions of $\frac{v_{1}(tN)}{\sqrt{N}}$ and $\frac{v_{2}(tN)}{\sqrt{N}}$,
where $v_{1}(t)=N_{t}^{+}-qt(1-\alpha) $ and $v_{2}(t)=N_{t}^{-}-pt(1-\alpha)$. With this purpose, we use Theorem 2.1 of \cite{S.}, with a
slightly different definition for weakly positively associated increments given in \cite{S.2}, namely:
\begin{definition}
A process $\{v(t): t\geq{0}\}$ has \textit{weakly positive associated} increments if for all coordinatewise increasing functions
$f:\mathbb{R}\rightarrow{\mathbb{R}}$, $g:\mathbb{R}^{n}\rightarrow{\mathbb{R}}$
\begin{equation*}
E_{\mu_{\alpha}}[f(v(t+s)-v(s))g(v(s_{1}),..,v(s_{n}))]\geq{E_{\mu_{\alpha}}[f(v(t))]E_{\mu_{\alpha}}[g(v(s_{1}),..,v(s_{n}))]},
\end{equation*}
for all $s,t\geq{0}$ and $0\leq{s_{1}}<..<s_{n}=s$ (weakly negative associated in the sense of the reversed inequality).
\end{definition}

Following the same arguments as in Theorem 2 of \cite{K.} we note that the processes $N_{t}^{+}$ and $N_{t}^{-}$, have weakly positive
associated increments. In the sake of completeness, we give a sketch of the proof of this result for the process $N_{t}^{+}$.

Let $s,t\geq{0}$ and $0\leq{s_{1}}<...<s_{n}=s$, and $f,g$ coordinatewise increasing functions $f:\mathbb{R}\rightarrow{\mathbb{R}}$,
$g:\mathbb{R}^{n}\rightarrow{\mathbb{R}}$. We have to show that
\begin{equation*}
E_{\mu_{\alpha}}[f(N_{t+s}^{+}-N_{s}^{+})g(N_{s_{1}}^{+},..,N_{s_{n}}^{+})]\geq{E_{\mu_{\alpha}}(f(N_{t}^{+}))E_{\mu_{\alpha}}(g(N_{s_{1}}^{+},..,N_{s_{n}}^{+}))}.
\end{equation*}
 Using the Markov property and by reverting the process $\eta_{s}$ with respect to $\mu_{\alpha}$ into $\hat{\eta_{s}}$, we have
\begin{equation*}
E_{\mu_{\alpha}}[f(N_{t+s}^{+}-N_{s}^{+})g(N_{s_{1}}^{+},..,N_{s_{n}}^{+})]=\int{E_{\eta}(f(N_{t}^{+}))\hat{E}_{\eta}(g(N_{s_{1}}^{-},..,N_{s_{n}}^{-}))d\mu_{\alpha}}.\label{eq:wpa}
\end{equation*}
Denote by $\varphi(\eta)$, $\psi(\eta)$ the functions $E_{\eta}(f(N_{t}^{+}))$ and $\hat{E}_{\eta}(g(N_{s_{1}}^{-},..,N_{s_{n}}^{-}))$,
respectively. Each one of this functions is increasing in each coordinate $\eta(x)$, because if we add one particle at site $x$, it can only
increase the number of jumps from $-1$ to $0$ (or from $0$ to $-1$). Using Lemma 3 of \cite{K.}, the right hand side of last expression is
bigger than
\begin{equation*}
\int{E_{\eta}(f(N_{t}^{+}))d\mu_{\alpha}}\int{\hat{E}_{\eta}(g(N_{s_{1}}^{-},...,N_{s_{n}}^{-}))d\mu_{\alpha}}.
\end{equation*}
And reversing the process again we obtain the result. For $N_{t}^{-}$ we can use the same argument.

Moreover, both processes have zero-mean and satisfy
\begin{equation*}
\lim_{t\rightarrow{+\infty}}\frac{1}{t}E_{\mu_{\alpha}}[(v_{i}(t))^{2}]=\sigma_{i}^{2}
\end{equation*}
for $i=1,2$ with $\sigma_{i}^{2}<\infty$, see Theorem 3 of \cite{K.}. In particular, the distributions of the processes
$\frac{v_{1}(tN)}{\sqrt{N}}$ and $\frac{v_{2}(tN)}{\sqrt{N}}$ are tight. The proof, see \cite{S.} relies on a maximal inequality, Corollary $6$
of \cite{N.W.}, which applies to demimartingales. As the processes have weakly positive associated increments and zero-mean, the demimartingale
property follows. \cqd

\section{Dependence on the initial configuration}

The first result we state concerns the dependence of the current through a fixed bond on the initial configuration. Here we suppose that $v>0$,
but for the other case, a similar statement holds.
\begin{proposition} \label{current dependence}
Fix $t\geq{0}$ and a site $x$. Then,
\begin{equation*}
\lim_{N\rightarrow{+\infty}}\mathbb{E}_{\nu_{\alpha}}\Big[\frac{\bar{J}^N_{x-1,x}(t)}{\sqrt{N}}-
\frac{\sum_{y=x-vtN}^{x-1}\bar{\eta}_{0}(y)}{\sqrt{N}}\Big]^{2}=0.
\end{equation*}
\end{proposition}
In the case $\alpha=1/2$, the normalized current converges to $0$ in the $L^{2}(\mathbb{P}_{\nu_{\alpha}})$-norm. This result was also obtained
before by Ferrari and Fontes in \cite{F.F.1}.
\begin{proof}
Here we consider $x=0$, but the same argument applied to any site $x$ provides the corresponding result. Recall the result of Proposition
\ref{prop:1}.

On the other hand, $Y^{N}_{t}(T_{t}G_{n})-Y^{N}_{0}(T_{t}G_{n})$ converges to $0$ as $N\rightarrow{+\infty}$ in the
$L^{2}(\mathbb{P}_{\nu_{\alpha}})$-norm, where $T_{t}H(u)=H(u+vt)$. For that, fix $H\in{S(\mathbb{R})}$, associate the martingales $M_{t}^{N,H}$
to the density fluctuation field, use the fact that $\mathbb{E}_{\nu_{\alpha}}[(M^{N,H}_{t})^2]$ vanishes as $N\rightarrow{+\infty}$ and the
Boltzmann-Gibbs Principle, see Theorem \ref{th:bg}. The result is accomplished for $G_{n}$, by approximating them in the
$L^{2}(\mathbb{P}_{\nu_{\alpha}})$-norm by smooth functions $H_{n,k}$ with compact support, as in the proof of Theorem 2.3 of \cite{J.L.}.

In order to finish the proof it remains to show that
\begin{equation*}
\lim_{n\rightarrow{+\infty}}\mathbb{E}_{\nu_{\alpha}}\Big[Y_{0}^{N}(T_{t}G_{n})-Y_{0}^{N}(G_{n})-
\frac{1}{\sqrt{N}}\sum_{x=-vtN}^{-1}\bar{\eta}_{0}(x)\Big]^{2}=0,
\end{equation*}
uniformly over $N$, which is a consequence of the explicit knowledge of $G_{n}$ and of $\nu_{\alpha}$ being a product measure.
\end{proof}
Since both, the current over a bond and the density fluctuation field at time $t$, can be written in terms of the initial configuration, and
since (\ref{eq:taggedparticle}) holds, it is natural that the position of the Tagged Particle also enjoys this property. That is the content of
the Corollary \ref{tagged dependence}, whose proof we start to present.

\quad\

\demon \\
We are going to show the convergence in $\mathbb{P}_{\nu_{\alpha}^*}$-probability to $0$ of the random variable appearing in the statement of
the Corollary, and then we show that its $L^2(\mathbb{P}_{\nu^*_{\alpha}})$-norm is finite, which allows to conclude the convergence to $0$ in
$L^{2-\epsilon}(\mathbb{P}_{\nu^*_{\alpha}})$, for any $\epsilon>0$.

With that purpose, start by summing and subtracting the expectation of $X_{tN}$, namely $v_{t}N$, in the expression that appears in the
statement of the Corollary, and it becomes as:
\begin{equation*}
\frac{\bar{X}_{tN}}{\sqrt{N}}+\frac{\sum_{x=1}^{(p-q)\alpha tN}\bar{\eta}_{0}(x)}{\alpha\sqrt{N}}.
\end{equation*}
We start by showing that last expression converges to zero in $\mathbb{P}_{\nu_{\alpha}^{*}}$-probability as $N\rightarrow{+\infty}$.

At first note that by the rigid transport of the system it holds that:
\begin{equation*}
\lim_{N\rightarrow{+\infty}}\mathbb{E}_{\nu_{\alpha}}\Big[\frac{\sum_{x=1+vtN}^{v_{t}N}\bar{\eta}_{t}(x)}{\alpha\sqrt{N}}-\frac{\sum_{x=1}^{(p-q)\alpha
tN}\bar{\eta}_{0}(x)}{\alpha\sqrt{N}}\Big]^2=0.
\end{equation*}
As a consequence we have to show that:
\begin{equation*}
\frac{\bar{X}_{tN}}{\sqrt{N}}+\frac{\sum_{x=1+vtN}^{v_{t}N}\bar{\eta}_{t}(x)}{\alpha\sqrt{N}},
\end{equation*}
converges to zero in $\mathbb{P}_{\nu_{\alpha}^{*}}$-probability as $N\rightarrow{+\infty}$.

In order to keep notation simple we denote by $Z^{N}_t$ the random variable:
\begin{equation*}
Z^{N}_t=-\sum_{x=1+vtN}^{v_{t}N}\bar{\eta}_{t}(x)/\alpha.
\end{equation*}

Notice that $v_{t}N+Z^{N}_t$ is a positive random variable since it corresponds to the number of holes in the interval $[1+vtN,v_{t}N]$.

Fix $a>0$, take $n=a\sqrt{N}+v_{t}N+Z^{N}_t$ in the expression that relates the position of the Tagged Particle with the current through the
bond $[-1,0]$ and the density of particles, see (\ref{eq:taggedparticle}):
\begin{equation*}
\{X_{tN}\geq{v_{t}N+a\sqrt{N}+Z^{N}_t}\}=
\Big\{J^N_{-1,0}(t)\geq{\sum_{x=0}^{v_{t}N}\eta_{t}(x)+\sum_{x=1+v_{t}N}^{a\sqrt{N}-1+v_{t}N+Z^{N}_t}\eta_{t}(x)}\Big\}.
\end{equation*}

Introducing the mean of the current, last expression becomes as
\begin{equation*}
\{X_{tN}\geq{v_{t}N+a\sqrt{N}+Z^{N}_t}\}=
\Big\{\bar{J}^N_{-1,0}(t)\geq{\sum_{x=0}^{v_{t}N}\bar{\eta}_{t}(x)+\sum_{x=1+v_{t}N}^{a\sqrt{N}-1+v_{t}N+Z^{N}_t}\eta_{t}(x)}\Big\}.
\end{equation*}
Now, we can divide all the terms by $\sqrt{N}$ and then, subtract the mean of the random variable on the right hand side of last inequality to
obtain:
\begin{equation*}
\{\frac{\bar{X}_{tN}}{\sqrt{N}}-\frac{Z^{N}_t}{\sqrt{N}}\geq{a}\}=
\end{equation*}
\begin{equation*}
\Big\{\frac{\bar{J}^N_{-1,0}(t)}{\sqrt{N}}\geq{\frac{\sum_{x=0}^{v_{t}N}\bar{\eta}_{t}(x)}{\sqrt{N}}+\frac{\sum_{x=1+v_{t}N}^
{a\sqrt{N}-1+v_{t}N+Z^{N}_t}\bar{\eta}_{t}(x)}{\sqrt{N}}+\alpha a +\frac{\alpha Z^{N}_t}{\sqrt{N}}}\Big\}.
\end{equation*}
By Proposition \ref{current dependence}, $T^N_{t}$ converges to zero in $L^2(\mathbb{P}_{\nu_{\alpha}})$ where
\begin{equation*}
T^{N}_{t}=\frac{\bar{J}^N_{-1,0}(t)}{\sqrt{N}}-\frac{1}{\sqrt{N}}\sum_{x=-vtN}^{-1}\bar{\eta}_{0}(x),
\end{equation*}
which together with the Boltzmann-Gibbs Principle gives us that:
\begin{equation*}
\mathbb{P}_{\nu_{\alpha}^{*}}\Big\{\frac{\bar{X}_{tN}}{\sqrt{N}}+\frac{\sum_{x=1+vtN}^{v_{t}N}\bar{\eta}_{t}(x)}{\alpha\sqrt{N}}\geq{a}\Big\}=
\end{equation*}
\begin{equation*}
\mathbb{P}_{\nu_{\alpha}^{*}}\Big\{\frac{\sum_{x=0}^{-1+vtN}\bar{\eta}_{t}(x)}{\sqrt{N}}\geq{\frac{\sum_{x=0}^{v_{t}N}\bar{\eta}_{t}(x)}{\sqrt{N}}+
\frac{\sum_{x=1+v_{t}N}^ {a\sqrt{N}-1+v_{t}N+Z^{N}_t}\bar{\eta}_{t}(x)}{\sqrt{N}}+\alpha a +\frac{\alpha Z^{N}_t}{\sqrt{N}}}\Big\}
\end{equation*}
Now observe that:
\begin{equation*}
\mathbb{E}_{\nu_{\alpha}}\Big[\frac{1}{\sqrt{N}}\sum_{x=1+v_{t}N}^ {a\sqrt{N}-1+v_{t}N+Z^{N}_t}\bar{\eta}_{t}(x)\Big]^2=O(N^{-1/2}),
\end{equation*}
whose proof is presented at the end in order to simplify the exposition. Therefore, for $N$ sufficiently big we have that
\begin{equation*}
\mathbb{P}_{\nu_{\alpha}^{*}}\{\frac{\bar{X}_{tN}}{\sqrt{N}}+\frac{\sum_{x=1+vtN}^{v_{t}N}\bar{\eta}_{t}(x)}{\alpha\sqrt{N}}\geq{a}\}=\mathbb{P}_{\nu_{\alpha}^{*}}\Big\{0
\geq{\frac{\sum_{x=vtN}^{v_{t}N}\bar{\eta}_{t}(x)}{\sqrt{N}}+ \alpha a -\frac{\sum_{x=1+vtN}^{v_{t} N}\bar{\eta}_{t}(x)}{\sqrt{N}}}\Big\}
\end{equation*}
which concludes the first step of the proof.

For the $L^{2-\epsilon}(\mathbb{P}_{\nu^*_{\alpha}})$ convergence, it remains to show that:
\begin{equation*}
\sup_{N}\mathbb{E}_{\nu_{\alpha}^{*}}\Big[\frac{\bar{X}_{tN}}{\sqrt{N}}+\frac{\sum_{x=1}^{(p-q)\alpha
tN}\bar{\eta}_{0}(x)}{\alpha\sqrt{N}}\Big]^{2}<{+\infty}.
\end{equation*}
Last result is a consequence of $\nu_{\alpha}$ being a product measure, which implies that:
\begin{equation*}
\mathbb{E}_{\nu_{\alpha}^{*}}\Big[\frac{\sum_{x=1}^{(p-q)\alpha tN}\eta_{0}(x)}{\sqrt{N}}\Big]^{2}\leq{(p-q)(1-\alpha) t};
\end{equation*}
together with a result due to De Masi and Ferrari in \cite{M.F.}:
\begin{equation*}
\lim_{N\rightarrow{+\infty}}\mathbb{E_{\nu_{\alpha}^{*}}}\Big[\frac{\bar{X}_{tN}}{\sqrt{N}}\Big]^{2}=(p-q)(1-\alpha)t.
\end{equation*}
In order to finish the proof it is enough to show that:
\begin{equation*}
\mathbb{E}_{\nu_{\alpha}}\Big[\frac{1}{\sqrt{N}}\sum_{x=1+v_{t}N}^ {a\sqrt{N}-1+v_{t}N+Z^{N}_t}\bar{\eta}_{t}(x)\Big]^2=O(N^{-1/2}).
\end{equation*}
To simplify the computations we take $p=1$, nevertheless the case $p\neq{1}$ follows the same lines. Since $\nu_{\alpha}$ is an invariant
measure, last expectation can be written as:
\begin{equation}
\int\frac{1}{N}\sum_{x,y=1+v_{t}N}^{a\sqrt{N}-1+v_{t}N+Z^{N}}\bar{\eta}(x)\bar{\eta}(y)\nu_{\alpha}(d\eta)=0, \label{eq1}
\end{equation}
where $Z^{N}$ is equal to:
\begin{equation*}
Z^{N}=-\sum_{x=1+vtN}^{v_{t}N}\bar{\eta}(x)/\alpha.
\end{equation*}
Notice that $Z^{N}$ depends on the variables $\eta(x)$ for $x$ depending on the sites from $1+vtN$ to $v_{t}N$, while the sum depends on the
random variables $\eta(x)$ for $x$ runing through $1+v_{t}N$ to $a\sqrt{N}-1+v_{t}N+Z^{N}$. So, we can separate the sum in (\ref{eq1}) into the
sites where the random variables appearing in the sum and $Z^{N}$ are independent from the sites where they correlate and it becomes as:
\begin{equation*}
\int_{\{a\sqrt{N}-1+v_{t}N+Z^{N}\geq{1+v_{t}N}\}}\frac{1}{N}\sum_{x,y=1+v_{t}N}^{a\sqrt{N}-1+v_{t}N+Z^{N}}\bar{\eta}(x)\bar{\eta}(y)\nu_{\alpha}(d\eta)=0
\end{equation*}
\begin{equation*}
+\int_{\{a\sqrt{N}-1+v_{t}N+Z^{N}<1+v_{t}N\}}\frac{1}{N}\sum_{x,y=1+v_{t}N}^{a\sqrt{N}-1+v_{t}N+Z^{N}}\bar{\eta}(x)\bar{\eta}(y)\nu_{\alpha}(d\eta)=0.
\end{equation*}
By independence the first integral is non zero as long as $x=y$ and it equals:
\begin{equation}
\frac{\alpha^2}{N}\int_{\{a\sqrt{N}+Z^{N}\geq{2}\}}\Big({a\sqrt{N}+Z^{N}-1}\Big)\nu_{\alpha}(d\eta), \label{eq:1}
\end{equation}
while the second can be written as:
\begin{equation}
\frac{1}{N}\int_{\{a\sqrt{N}+Z^{N}<2\}}\sum_{x,y=a\sqrt{N}-1+v_{t}N+Z^{N}}^{1+v_{t}N}\eta(x)\eta(y)\nu_{\alpha}(d\eta) \label{int:1}
\end{equation}
\begin{equation}
-\frac{2\alpha}{N}\int_{\{a\sqrt{N}+Z^{N}<2\}}\sum_{x,y=a\sqrt{N}-1+v_{t}N+Z^{N}}^{1+v_{t}N}\eta(x)\nu_{\alpha}(d\eta) \label{int:2}
\end{equation}
\begin{equation}
+\frac{\alpha^2}{N}\int_{\{a\sqrt{N}+Z^{N}<2\}}\sum_{x,y=a\sqrt{N}-1+v_{t}N+Z^{N}}^{1+v_{t}N}\nu_{\alpha}(d\eta). \label{int:3}
\end{equation}
Now, we give the route to proceed in the computations. For $j=1,2$, let $Z^{N,j}$ be the random variable:
\begin{equation*}
Z^{N,j}=-\sum_{x=1+vtN}^{v_{t}N-j}\bar{\eta}(x)/\alpha.
\end{equation*}
Estimate (\ref{int:1}) by separating the case $x=y$ from the case $x\neq{y}$. In the first one the integral becomes as:
\begin{equation*}
\frac{\alpha}{N}\int_{\{a\sqrt{N}+Z^{N,1}<2+(1-\alpha)/\alpha\}}\Big(2+\frac{(1-\alpha)}{\alpha}-a\sqrt{N}-Z^{N,1}\Big)\nu_{\alpha}(d\eta),
\end{equation*}
while in the case $x\neq{y}$ it becomes as:
\begin{equation*}
\frac{\alpha^2}{N}\int_{\{a\sqrt{N}+Z^{N,2}<2+2(1-\alpha)/\alpha\}}\Big(2+\frac{2(1-\alpha)}{\alpha}-a\sqrt{N}-Z^{N,2}\Big)\Big(3-a\sqrt{N}-Z^{N,2}\Big)\nu_{\alpha}(d\eta).
\end{equation*}
On the other hand, (\ref{int:2}) can be written as
\begin{equation*}
\frac{-\alpha}{N}\int_{\{a\sqrt{N}+Z^{N,1}<2+(1-\alpha)/\alpha\}}\Big(2+\frac{(1-\alpha)}{\alpha}-a\sqrt{N}-Z^{N,1}\Big)^2\nu_{\alpha}(d\eta),
\label{eq:4}
\end{equation*}
while (\ref{int:3}) is equal to
\begin{equation}
\frac{\alpha^2}{N}\int_{\{a\sqrt{N}+Z^{N}<2\}}\Big(2-a\sqrt{N}-Z^{N}\Big)^2\nu_{\alpha}(d\eta). \label{eq:5}
\end{equation}
Now, it remains to write all the integrals with respect to the random variable $Z^{N,2}$. Since the Bernoulli product measure is homogenous we
condition on $\eta(x)=0$ and $\eta(x)=1$ for some site $x\in[1+vtN,v_{t}N]$, to write the integrals (\ref{eq:1}) and (\ref{eq:5}) in terms of
$Z^{N,1}$. Them we repeat the same procedure to write the remaining integrals in terms of $Z^{N,2}$. Organizing them all, the result follows.
 \cqd
\section{Density Fluctuations in a longer time scale}

Here we are focused in proving Theorem \ref{th:flu2}. Fix a positive integer $k$ and recall the definition of the density fluctuation field in
(\ref{eq:densfieldlongscale}). Let $U_{t}^{N}H(u)=H(u-vtN^{\gamma})$. As before, we need to show that $(Q^\gamma_{N})_{N\geq{1}}$ is tight and
to characterize the limit field. We start by the latter while the former is referred to the eighth section.

Fix $H\in{S(\mathbb{R})}$. Then:
\begin{equation} \label{martingaleM}
M_{t}^{N,H}=Y_{t}^{N,\gamma}(H)-Y_{0}^{N,\gamma}(H)-\int_{0}^{t}\Gamma^{H}_{1}(s)ds,
\end{equation}
is a martingale with respect to $\mathcal{\tilde{F}}_{t}=\sigma(\eta_{s}, s\leq{t})$ whose quadratic variation is given by
\begin{equation} \label{quadraticlonger}
\int^{t}_{0}\frac{N^{\gamma}}{N^2}\sum_{x\in\mathbb{Z}}
\Big(\nabla^{N}{U_{s}^{N}H\Big(\frac{x}{N}\Big)}\Big)^2\Big[c(x,x+1,\eta_{s})+c(x+1,x,\eta_{s})\Big]ds,
\end{equation}
where $\Gamma^{H}_{1}(s)$ equals to
\begin{equation} \label{intpartofmartingale}
\frac{N^{\gamma}}{\sqrt{N}}\sum_{x\in{\mathbb{Z}}}\nabla^{N} U_{s}^{N}H\Big(\frac{x}{N}\Big)W_{x,x+1}(\eta_{s})-
\frac{N^{\gamma}}{\sqrt{N}}\sum_{x\in{\mathbb{Z}}}\partial_{u}U_{s}^{N}H\Big(\frac{x}{N}\Big)v[\eta_{s}(x)-\alpha].
\end{equation}

Easily one shows that the $L^{2}(\mathbb{P}_{\nu_\alpha}^{\gamma})$-norm of $M_{t}^{N,H}$ vanishes as $N\rightarrow{+\infty}$ as long as
$\gamma<1$. Then, under a sub-diffusive time scale regime, the only term contributing to the limit density fluctuation field is its integral
part, since its quadratic variation vanishes. The characterization of the limit of the integral part of the martingale is known as the
Boltzmann-Gibbs Principle and is the main difficulty when showing the equilibrium fluctuations. In that scaling regime the time evolution of the
limit density fluctuation field is given in a similar way to (\ref{timeevolutionofdfield}). But when one takes the diffusive scaling a new
contribution arises, since the quadratic variation of the martingale does not vanishes, which agrees with the fact that in order to observe
fluctuations from the dynamics one has to take this time scale.

Now, we proceed by proving that the integral part of the martingale $M_{t}^{N,H}$ vanishes in $L^{2}(\mathbb{P}_{\nu_\alpha}^{\gamma})$ as
$N\rightarrow{+\infty}$. Since $\sum_{x\in{\mathbb{Z}}}\nabla^N U_{s}^{N}H\Big(\frac{x}{N}\Big)=0$, we can introduce it times
$E_{\nu_{\alpha}}[W_{x,x+1}(\eta)]$, in the integral part of the martingale $M_{t}^{N,H}$ and using the decomposition of the instantaneous
current
\begin{equation} \label{eq:decompcurrent}
\bar{W}_{0,1}(\eta)=-(p-q)\bar{\eta}(0)\bar{\eta}(1)-(q(1-\alpha)+p\alpha)[\bar{\eta}(1)-\bar{\eta}(0)]+v[\eta(0)-\alpha],
\end{equation}
it becomes as:
\begin{equation*}
\int^{t}_{0}\frac{N^{\gamma}}{\sqrt{N}}\sum_{x\in\mathbb{Z}}\nabla^{N}U_{s}^{N}
H\Big(\frac{x}{N}\Big)(q-p)\bar{\eta}_{s}(x)\bar{\eta}_{s}(x+1)ds
\end{equation*}
\begin{equation*}
+\int^{t}_{0}\frac{N^{\gamma}}{\sqrt{N}}\sum_{x\in\mathbb{Z}}\nabla^{N}
U_{s}^{N}H\Big(\frac{x}{N}\Big)(q(1-\alpha)+p\alpha)[\bar{\eta}_{s}(x+1)-\bar{\eta}_{s}(x)]ds
\end{equation*}
\begin{equation*}
+\int^{t}_{0}\frac{N^{\gamma}}{\sqrt{N}}\sum_{x\in\mathbb{Z}}\Big\{\nabla^{N}U_{s}^{N}H\Big(\frac{x}{N}\Big)-
\partial_{u}U_{s}^{N}H\Big(\frac{x}{N}\Big)\Big\}v[\eta_{s}(x)-\alpha]ds.
\end{equation*}

By a summation by parts, Schwarz inequality and since $\nu_{\alpha}$ is a product invariant measure, the second term vanishes in
$L^{2}(\mathbb{P}_{\nu_{\alpha}}^{\gamma})$ as $N\rightarrow{+\infty}$, while for the last term we use Taylor expansion to show that it vanishes
in the same norm. Once more, last results hold as long as $\gamma<1$.

It remains to show that the $L^{2}(\mathbb{P}_{\nu_\alpha}^{\gamma})$-norm of the first integral vanishes as $N\rightarrow{+\infty}$. For that,
we use the Botzmann-Gibbs Principle, which is proved in the next section. This result is accomplished for $\gamma<1/3$, but it should hold for
$\gamma<1/2$ as conjectured. We also remark, that almost all the subsequent results rely on the Boltzmann-Gibbs Principle and if one shows that
it holds for $\gamma<1/2$, one can establish the same results up to the time scale $N^{3/2}$.

Assuming that $(Q^\gamma_{N})_{N}$ is tight, it has convergent subsequences.
 Let $Q$ be one of its limiting points. By the results proved so far, under
$Q$, the density fluctuation field satisfies $Y_{t}(H)=Y_{0}(H)$.

 For $t\geq0$, let $\mathcal{F}_{t}$ be the $\sigma$-algebra on $D([0,T],\mathcal{H}_{-k})$ generated by $Y_{s}(H)$ for
$s\leq{t}$ and $H$ in $S(\mathbb{R})$ and set $\mathcal{F}=\sigma(\bigcup_{t\geq{0}}\mathcal{F}_{t})$. It is not hard to show
 as in chap. 11 of \cite{K.L.}, that up to this longer
time scale $N^{4/3}$, $Q$ restricted to $\mathcal{F}_{0}$ is a Gaussian field with covariance given by
$E_{Q}(Y_{0}(G)Y_{0}(H))=\chi(\alpha)<G,H>$ and it is trivial that the limit density field has covariance given by (\ref{eq:covariancelonger}).
This concludes the proof of the Theorem \ref{th:flu2}.

\section{Boltzmann-Gibbs Principle}

In this section we prove Theorem \ref{bg}.

Fix $H\in S(\mathbb{\mathbb{R}})$ and an integer $K$. We divide $\mathbb{Z}$ in non overlapping intervals of length $K$, denoted by
$\{I_{j},j\geq{1}\}$. Then, the expectation can be written as:
\begin{equation*}
\mathbb{E}_{\nu_{\alpha}}^{\gamma}\Big[\int_{0}^{t}\frac{N^{\gamma}}{\sqrt{N}}\sum_{j\geq{1}}\sum_{x\in{I}_{j}}H\Big(\frac{x}{N}\Big)\bar{\eta}_{s}(x)\bar{\eta}_{s}(x+1)ds\Big]^{2}.
\end{equation*}

In order to have independence of $\bar{\eta}(x)\bar{\eta}(x+1)$ and $\bar{\eta}(y)\bar{\eta}(y+1)$ for $x$ and $y$ in different $I_{j}$'s, we
separate the sum over the intervals $I_{j}$ for $j$ odd, and $j$ even. So, in fact it remains to bound
\begin{equation} \label{eq:bg1}
\mathbb{E}_{\nu_{\alpha}}^{\gamma}\Big[\int_{0}^{t}\frac{N^{\gamma}}{\sqrt{N}}\sum_{j\geq{1}}\sum_{x\in{I}_{j}}
H\Big(\frac{x}{N}\Big)\bar{\eta}_{s}(x)\bar{\eta}_{s}(x+1)ds\Big]^{2}
\end{equation}
where, for example $j$ is odd. The case for $j$ even follows by the same arguments. Remark that, in this setting, every $x\in{I_{j}}$ and
$y\in{I_{l}}$, for $j\neq{l}$, are at least at a distance $K$.

Now, sum and subtract $H\Big(\frac{y_{j}}{N}\Big)$, where $y_{j}$ is a point of the interval $I_{j}$, inside the summation over $x$. Since
$(x+y)^2\leq{2x^2+2y^2}$, the expression (\ref{eq:bg1}) can be bounded by

\begin{equation*}
2\mathbb{E}_{\nu_{\alpha}}^{\gamma}\Big[\int_{0}^{t}\frac{N^{\gamma}}{\sqrt{N}}\sum_{j\geq{1}}
\sum_{x\in{I}_{j}}\Big[H\Big(\frac{x}{N}\Big)-H\Big(\frac{y_{j}}{N}\Big)\Big]\bar{\eta}_{s}(x)\bar{\eta}_{s}(x+1)ds\Big]^{2}
\end{equation*}
\begin{equation}\label{eq:bg3}
+2\mathbb{E}_{\nu_{\alpha}}^{\gamma}\Big[\int_{0}^{t}\frac{N^{\gamma}}{\sqrt{N}}\sum_{j\geq{1}}
H\Big(\frac{y_{j}}{N}\Big)\sum_{x\in{I_{j}}}\bar{\eta}_{s}(x)\bar{\eta}_{s}(x+1)ds\Big]^{2}.
\end{equation}
We are going to estimate each term separately and divide the proof in several lemmas, to make the exposition clearer. We start by the former.
\begin{lemma} \label{lemma:lenghK1}
For every $H\in S(\mathbb{R})$ and every $t>0$, if $KN^{\gamma-1}\rightarrow{0}$ as $N\rightarrow{+\infty}$, then
\begin{equation*}
\lim_{N\rightarrow{\infty}}\mathbb{E}_{\nu_{\alpha}}^{\gamma}\Big[\int_{0}^{t}\frac{N^{\gamma}}{\sqrt{N}}
\sum_{j\geq{1}}\sum_{x\in{I}_{j}}\Big[H\Big(\frac{x}{N}\Big)-H\Big(\frac{y_{j}}{N}\Big)\Big]\bar{\eta}_{s}(x)\bar{\eta}_{s}(x+1)ds\Big]^{2}=0.
\end{equation*}
\end{lemma}
\begin{proof}
By Schwarz inequality and since $\nu_{\alpha}$ is an invariant product measure, the expectation is bounded by $C
t^{2}\frac{N^{2\gamma}}{N}\sum_{j}\sum_{x\in{I_{j}}}\Big(H'\Big(\frac{y_{j}}{N}\Big)\Big)^{2}\Big(\frac{|x-y_{j}|}{N}\Big)^{2}$. Since $x$ and
$y_{j}$ are in the $I_{j}$ interval, that has size $K$, last expression can be bounded by $Ct^{2}N^{2\gamma}||H'||_{2}^{2}(\frac{K}{N})^{2}$
which vanishes as long as $K N^{\gamma-1}\rightarrow{0}$  when $N\rightarrow{+\infty}$.
\end{proof}
 Now, we bound the expression (\ref{eq:bg3}). We sum and subtract the expectation of $\sum_{x\in{I_{j}}}\bar{\eta}_{s}(x)\bar{\eta}_{s}(x+1)$
conditioned on the hyperplanes $M_{j}=\sigma\Big(\sum_{x\in{I^{*}_{j}}}\eta(x)\Big)$, where $I_{j}^{*}=I_{j}\bigcup\{x_{j+1}\}$, if
$I_{j}=\{x_{0},x_{1},..,x_{j}\}$. Using again the elementary inequality $(x+y)^2\leq{2x^2+2y^2}$, the expectation in (\ref{eq:bg3}) is bounded
by
\begin{equation} \label{eq:bg4}
2\mathbb{E}_{\nu_{\alpha}}^{\gamma}\Big[\int_{0}^{t}\frac{N^{\gamma}}{\sqrt{N}}\sum_{j\geq{1}}H\Big(\frac{y_{j}}{N}\Big)V_{j}(\eta_{s})ds\Big]^{2}
\end{equation}
\begin{equation}\label{eq:bg5}
+2\mathbb{E}_{\nu_{\alpha}}^{\gamma}\Big[\int_{0}^{t}\frac{N^{\gamma}}{\sqrt{N}}\sum_{j\geq{1}}H\Big(\frac{y_{j}}{N}\Big)E\Big(\sum_{x\in{I_{j}}}
\bar{\eta}_{s}(x)\bar{\eta}_{s}(x+1)\Big|M_{j}\Big)ds\Big]^{2}
\end{equation}
where
\begin{equation*}
V_{j}(\eta)=\sum_{x\in{I_{j}}} \bar{\eta}(x)\bar{\eta}(x+1)-E\Big(\sum_{x\in{I_{j}}} \bar{\eta}(x)\bar{\eta}(x+1)\Big|M_{j}\Big).
\end{equation*}
Once more, we bound the integrals separately. We start by bounding (\ref{eq:bg4}).
\begin{lemma} \label{lm:lengthK}
For every $H\in S(\mathbb{R})$ and every $t>0$, if $K^{2}N^{\gamma-1}\rightarrow{0}$ as $N\rightarrow{+\infty}$, then
\begin{equation*}
\lim_{N\rightarrow{\infty}}\mathbb{E}_{\nu_{\alpha}}^{\gamma}\Big[\int_{0}^{t}\frac{N^{\gamma}}{\sqrt{N}}\sum_{j\geq{1}}H\Big(\frac{y_{j}}{N}\Big)
V_{j}(\eta_{s})ds\Big]^{2}=0.
\end{equation*}
\end{lemma}
\begin{proof}
For $f,g\in{L^{2}(\nu_{\alpha})}$ define the inner product $<f,-Lg>_{\nu_{\alpha}}$.  Let ${H}_{1}$ be the Hilbert space generated by
$L^{2}(\nu_{\alpha})$ and this inner product. Denote by $||\cdot||_{1}$ the norm induced by this inner product and let $||\cdot||_{-1}$ be its
dual norm with respect to $L^{2}(\nu_{\alpha})$:
\begin{equation} \label{variationalformH-1norm}
||f||_{-1}=\sup_{g\in{L^{2}(\nu_{\alpha})}}\Big\{2<f,g>_{\nu_{\alpha}}-||g||_{1}\Big\}.
\end{equation}
By definition for every $f\in{H_{-1}}$, $g\in{L^{2}(\nu_{\alpha})}$ and $A>0$ it holds that:
\begin{equation} \label{norm bound}
2<f,g>_{\nu_{\alpha}}\leq{\frac{1}{A}||f||_{-1}+A||g||_{1}}.
\end{equation}
By Proposition A1.6.1 of \cite{K.L.}, the expectation in the statement of the Lemma is bounded by
\begin{equation*}
Ct\Big|\Big|\frac{N^{\gamma}}{\sqrt{N}}\sum_{j\geq{1}}H\Big(\frac{y_{j}}{N}\Big)V_{j}\Big|\Big|_{-1}^{2},
\end{equation*}
where $C$ is a constant. By the variational formula for the $H_{-1}$-norm (\ref{variationalformH-1norm}) last expression is equal to
\begin{equation*}
Ct \sup_{h \in L^{2}(\nu_{\alpha})}\Big\{2\int
\frac{N^{\gamma}}{\sqrt{N}}\sum_{j\geq{1}}H\Big(\frac{y_{j}}{N}\Big)V_{j}(\eta)h(\eta)\nu_{\alpha}(d\eta)-N^{1+\gamma}<h,-L^{S}_{N}h>_{\alpha}\Big\},
\end{equation*}
and is bounded by
\begin{equation*}
Ct \sum_{j\geq{1}}\sup_{h \in L^{2}(\nu_{\alpha})}\Big\{2\int
\frac{N^{\gamma}}{\sqrt{N}}H\Big(\frac{y_{j}}{N}\Big)V_{j}(\eta)h(\eta)\nu_{\alpha}(d\eta)-N^{1+\gamma}<h,-L^{S}_{I^{*}_{j}}h>_{\alpha}\Big\},
\end{equation*}
where $L^{S}_{I^{*}_{j}}$ denotes the restriction of the generator of the SSEP that we denote by $L^{S}_{N}$, to the set $I_{j}^{*}$:
\begin{equation*}
L^{S}_{I^{*}_{j}}f(\eta)=\sum_{\substack{x,y\in{I_{j}^{*}}\\|x-y|=1}}\frac{1}{2}\eta(x)(1-\eta(y))[f(\eta^{x,y})-f(\eta)].
\end{equation*}

Since $E(V_{j}|M_{j})=0$, $V_{j}$ belongs to the image of the generator $L^{S}_{I^{*}_{j}}$. Therefore, by (\ref{norm bound}) for each $j$ and
$A_{j}$ a positive constant it holds that
\begin{equation*}
\int V_{j}(\eta)h(\eta)\nu_{\alpha}(d\eta)\leq{\frac{1}{2
A_{j}}<V_{j},(-L^{S}_{I_{j}^{*}})^{-1}V_{j}>_{\alpha}+\frac{A_{j}}{2}<h,-L^{S}_{I_{j}^{*}}h>_{\alpha}}.
\end{equation*}
Taking for each $j$, $A_{j}=N^{3/2}\Big(|H(\frac{y_{j}}{N})|\Big)^{-1}$, the expectation becomes bounded by
\begin{equation*}
Ct\sum_{j\geq{1}}\frac{N^{\gamma}}{N^{2}}{H^{2}\Big(\frac{y_{j}}{N}\Big)}<V_{j},(-L^{S}_{I_{j}^{*}})^{-1}V_{j}>_{\alpha},
\end{equation*}
since the other term cancels with the $H_{1}$-norm of $h$. By the spectral gap inequality for the SSEP (see \cite{Q}) last expression can be
bounded by
\begin{equation*}
Ct\sum_{j\geq{1}}\frac{N^{\gamma}}{N^{2}}{H^{2}\Big(\frac{y_{j}}{N}\Big)}(K+1)^{2}Var(V_{j},\nu_{\alpha}).
\end{equation*}
 Now we observe that, since we are considering the extended interval $I_{j}^{*}$, it holds that
\begin{equation*}
E\Big(\bar{\eta}(x)\bar{\eta}(x+1)\Big|M_{j}\Big)=(\eta^{K+1}-\alpha)^{2}-\frac{1}{K}\eta^{K+1}(1-\eta^{K+1}),
\end{equation*}
where $\eta^{K+1}=(K+1)^{-1}\sum_{x\in{I_{j}^*}}\eta(x)$.

By a simple computation it is not hard to show that $Var(V_{j},\nu_{\alpha})\leq{K C}$, which implies the integral to be bounded by $C
t\frac{N^{\gamma}}{N}(K+1)^{2}||H||_{2}^{2}$ and vanishes as long as $K^{2}N^{\gamma-1}\rightarrow{0}$
 when $N\rightarrow{+\infty}$.
\end{proof}
To conclude the proof of the theorem it remains to bound (\ref{eq:bg5}). The idea we use to proceed consists in doing the following. Fix an
integer $L$ and consider bigger disjoint intervals of length $M=LK$, denoted by $\{\tilde{I}_{l},l\geq{1}\}$. In this setting, we consider $L$
sets of size $K$ together and we are able to write the expectation appearing in (\ref{eq:bg5}) as:
\begin{equation*}
\mathbb{E}_{\nu_{\alpha}}^{\gamma}\Big[\int_{0}^{t}\frac{N^{\gamma}}{\sqrt{N}}\sum_{l\geq{1}}\sum_{j\in{\tilde{I}_{l}}}H\Big(\frac{y_{j}}{N}\Big)
E\Big(\sum_{x\in{I_{j}}}\bar{\eta}_{s}(x)\bar{\eta}_{s}(x+1)\Big|M_{j}\Big)ds\Big]^{2}.
\end{equation*}
As before, sum and subtract $H\Big(\frac{z_{l}}{N}\Big)$, where $z_{l}$ denotes one point of the interval $\tilde{I}_{l}$, inside the summation
over $j$. Since $(x+y)^2\leq{2x^2+2y^2}$, last expectation can be bounded by
\begin{equation*}
2\mathbb{E}_{\nu_{\alpha}}^{\gamma}\Big[\int_{0}^{t}\frac{N^{\gamma}}{\sqrt{N}}\sum_{l\geq{1}}\sum_{j\in{\tilde{I}_{l}}}\Big[H\Big(\frac{y_{j}}{N}\Big)-H\Big(\frac{z_{l}}{N}\Big)\Big]
E\Big(\sum_{x\in{I_{j}}}\bar{\eta}_{s}(x)\bar{\eta}_{s}(x+1)\Big|M_{j}\Big)ds\Big]^{2}
\end{equation*}
\begin{equation} \label{eq:bg7}
+2\mathbb{E}_{\nu_{\alpha}}^{\gamma}\Big[\int_{0}^{t}\frac{N^{\gamma}}{\sqrt{N}}\sum_{l\geq{1}}H\Big(\frac{z_{l}}{N}\Big)\sum_{j\in{\tilde{I}_{l}}}E\Big(\sum_{x\in{I_{j}}}
\bar{\eta}_{s}(x)\bar{\eta}_{s}(x+1)\Big|M_{j}\Big)ds\Big]^{2}.
\end{equation}
The first expectation can be treated in the way as in the proof of Lemma \ref{lemma:lenghK1}, and it vanishes if
 $L^{2}KN^{2\gamma-2}\rightarrow{0}$ as $N\rightarrow{+\infty}$.

For the remaining expectation (\ref{eq:bg7}), inside the sum over $l$, sum and subtract
$E\Big(\sum_{x\in{\tilde{I}_{l}}}\bar{\eta}(x)\bar{\eta}(x+1)\Big|\tilde{M}_{l}\Big)$ where
$\tilde{M}_{l}=\sigma\Big(\sum_{x\in{\tilde{I}^{*}_{l}}}\eta(x)\Big)$ and $\tilde{I}_{l}^{*}$ denotes the extended interval $\tilde{I}_{l}$.
Then, the expectation in (\ref{eq:bg7}) can be bounded by
\begin{equation}\label{eq:bg8}
2\mathbb{E}_{\nu_{\alpha}}^{\gamma}\Big[\int_{0}^{t}\frac{N^{\gamma}}{\sqrt{N}}\sum_{l\geq{1}}H\Big(\frac{z_{l}}{N}\Big)\tilde{V}_{l}(\eta_{s})ds\Big]^{2}
\end{equation}
\begin{equation} \label{eq:bg9}
+2\mathbb{E}_{\nu_{\alpha}}^{\gamma}\Big[\int_{0}^{t}\frac{N^{\gamma}}{\sqrt{N}}\sum_{l\geq{1}}H\Big(\frac{z_{l}}{N}\Big)
E\Big(\sum_{x\in{\tilde{I}_{l}}}\bar{\eta}_{s}(x)\bar{\eta}_{s}(x+1)\Big|\tilde{M}_{l}\Big)ds\Big]^{2},
\end{equation}
 where
\begin{equation*}
\tilde{V}_{l}(\eta)=\sum_{j\in{\tilde{I}_{l}}}E\Big(\sum_{x\in{I_{j}}}
\bar{\eta}(x)\bar{\eta}(x+1)\Big|M_{j}\Big)-E\Big(\sum_{x\in{\tilde{I}_{l}}}\bar{\eta}(x)\bar{\eta}(x+1)\Big|\tilde{M}_{l}\Big).
\end{equation*}
We proceed by estimating (\ref{eq:bg8}):

\begin{lemma} \label{lm:lengthM}
For every $H\in S(\mathbb{R})$ and every $t>0$, if $L^{2}KN^{\gamma-1}\rightarrow{0}$ as $N\rightarrow{+\infty}$, then
\begin{equation*}
\lim_{N\rightarrow{\infty}}\mathbb{E}_{\nu_{\alpha}}^{\gamma}\Big[\int_{0}^{t}\frac{N^{\gamma}}{\sqrt{N}}\sum_{l\geq{1}}H\Big(\frac{z_{l}}{N}\Big)
\tilde{V}_{l}(\eta_{s})ds\Big]^{2}=0.
\end{equation*}
\end{lemma}
\begin{proof}
 Using the same arguments as in the proof of Lemma \ref{lm:lengthK}, the expectation becomes bounded by
\begin{equation*}
Ct \sum_{l\geq{1}}\sup_{h \in L^{2}(\nu_{\alpha})}\Big\{2 \int \frac{N^{\gamma}}{\sqrt{N}}
H\Big(\frac{z_{l}}{N}\Big)\tilde{V}_{l}(\eta)h(\eta)\nu_{\alpha}(d\eta)-N^{1+\gamma}<h,-L^{S}_{\tilde{I}_{l}^{*}}h>_{\alpha}\Big\}.
\end{equation*}
 Using an appropriate $A_{l}$ and the spectral gap inequality, we can bound last expression by
\begin{equation*}
Ct\sum_{l\geq{1}}\frac{N^{\gamma}}{N^{2}}H^{2}\Big(\frac{z_{l}}{N}\Big)(M+1)^{2}Var(\tilde{V}_{l},\nu_{\alpha}),
\end{equation*}
and since $Var(\tilde{V}_{l},\nu_{\alpha})\leq{L C}$ it vanishes if $L^{2}KN^{\gamma-1}\rightarrow{0}$, as $N\rightarrow{+\infty}$.
\end{proof}
To treat the remaining expectation (\ref{eq:bg9}) we continue applying the same steps.

\vspace{0,2cm}

\textbf{The proof of Boltzmann-Gibbs Principle}

\vspace{0,2cm} The idea of the proof was to take intervals of growing size in each step, in a way that the expectation vanishes for certain
restrictions on this size. The size of the first intervals taken, was $K$ and the biggest restriction in this size comes from Lemma
\ref{lm:lengthK}, namely that $K$ is such that $K^{2}N^{1-\gamma}\rightarrow{0}$ as $N\rightarrow{+\infty}$. Therefore, we can take
$K=N^{\frac{1-\gamma}{2}-\epsilon}$.

In the second step we had intervals of bigger size, namely $M$, where $M=LK$ and the parameter $L$ has to satisfy
$L^{2}KN^{\gamma-1}\rightarrow{0}$ as $N\rightarrow{+\infty}$. Since in the first step $K=N^{\frac{1-\gamma}{2}-\epsilon}$, we can take
$L=N^\frac{1-\gamma}{4}$, and as a consequence $M=N^{\frac{1-\gamma}{2}+{\frac{1-\gamma}{4}}-\epsilon}$.

Continuing the proof applying the same arguments, in the $n^{th}$ step we have intervals, denoted by $\{I^{n}_{p},p\geq{1\geq}\}$ of length
$K_{n}=N^{a_{n}}$, where $a_{n}=(1-\gamma)(\frac{1}{2}+\frac{1}{2^{2}}+...+\frac{1}{2^{n}})-\epsilon$.

Supposing that we stop this induction procedure in the $n^{th}$ step, it remains to bound the following expectation:
\begin{equation*}
\mathbb{E}_{\nu_{\alpha}}^{\gamma}\Big[\int_{0}^{t}\frac{N^{\gamma}}{\sqrt{N}}\sum_{p\geq{1}}H\Big(\frac{z_{p}}{N}\Big)
E\Big(\sum_{x\in{{I}^{n}_{p}}}\bar{\eta}_{s}(x)\bar{\eta}_{s}(x+1)\Big|{M}_{p}^{n}\Big)ds\Big]^{2},
\end{equation*}
where for each $p$, ${I}^{n}_{p}$ is an interval of size $K_{n}$, $z_{p}$ is one point of it and the hyperplanes are
${M}_{p}^{n}=\sigma\Big(\sum_{x\in{({I}^{n}_{p})^{*}}}\eta(x)\Big)$, where $({I}^{n}_{p})^{*}$ is taken as above.

Since $\nu_{\alpha}$ is an invariant product measure, last expectation can be bounded by
\begin{equation*}
t^{2}\frac{N^{2\gamma}}{N}\sum_{p\geq{1}}\Big(H\Big(\frac{z_{p}}{N}\Big)\Big)^{2}{E}_{\nu_{\alpha}}\Big(
E\Big(\sum_{x\in{{I}^{n}_{p}}}\bar{\eta}(x)\bar{\eta}(x+1)\Big|{M}_{p}^{n}\Big)\Big)^{2}.
\end{equation*}
Now, it is not hard to show that $E_{\nu_{\alpha}}\Big(
E\Big(\sum_{x\in{{I}^{n}_{p}}}\bar{\eta}(x)\bar{\eta}(x+1)\Big|M_{p}^{n}\Big)\Big)^{2}=O(1)$. Then the integral becomes bounded by
$\frac{N^{2\gamma}}{K_{n}}$, and for $n$ sufficiently big, since $K_{n}\sim{N^{1-\gamma}}$ and $\gamma<1/3$ this expression vanishes as
$N\rightarrow{+\infty}$. Here is the point in the proof where we need to impose the restriction on the parameter $\gamma<1/3$.
\begin{remark}
 Here we give an application of the Boltzmann-Gibbs Principle for the quadratic density fluctuation field associated to the
one-dimensional SSEP, in the diffusive scaling.

Consider a Markov process $\eta_{t}^s$ with generator given by (\ref{eq:generator}), with $p(x,y)=1/2$ under diffusive time scale. Consider
$\mathbb{P}^{N}_{\nu_{\alpha}}=\mathbb{P}_{\nu_{\alpha}}$ the probability measure on $D(\mathbb{R}^{+},\{0,1\}^{\mathbb{Z}})$ induced by the
invariant measure $\nu_{\alpha}$ and the Markov process $\eta_{t}^s$ speeded up by $N^2$ and denote by $\mathbb{E}_{\nu_{\alpha}}$ the
expectation with respect to $\mathbb{P}_{\nu_{\alpha}}$.

 Define the
quadratic density fluctuation field on $H\in{S(\mathbb{R})}$ by:
\begin{equation*}
\mathcal{Y}_{t}^{N}(H)=\frac{1}{\sqrt{N}}\sum_{x\in{\mathbb{Z}}}H\Big(\frac{x}{N}\Big)[\eta_{tN^2}^s(x)-\alpha][\eta_{tN^2}^s(x+1)-\alpha].
\end{equation*}
Following the same steps as in the proof of the Boltzmann-Gibbs Principle it is easy to show that:
\begin{corollary}
Fix $t>0$ and $\beta<1/2$, then
\begin{equation*}
\lim_{N\rightarrow{\infty}}\mathbb{E}_{\nu_{\alpha}}\Big[{N^{\beta}}\int_{0}^{t}
\frac{1}{\sqrt{N}}\sum_{x\in{\mathbb{Z}}}H\Big(\frac{x}{N}\Big)[\eta_{tN^2}^s(x)-\alpha][\eta_{tN^2}^s(x+1)-\alpha]ds\Big]^{2}=0.
\end{equation*}
\end{corollary}
Therefore, in order to observe fluctuations for the quadratic density fluctuation field, we need to consider $\beta\geq{1/2}$. In fact, in
\cite{Assing.} it is shown that
\begin{equation*}
{N^{1/2}}\int_{0}^{t} \frac{1}{\sqrt{N}}\sum_{x\in{\mathbb{Z}}}H\Big(\frac{x}{N}\Big)[\eta_{tN^2}^s(x)-\alpha][\eta_{tN^2}^s(x+1)-\alpha]ds
\end{equation*}
converges in law to a non-Gaussian singular functional of an infinite Ornstein-Uhlenbeck process.
\end{remark}
\section{Tightness}

 Now we prove that the sequence of probability measures $(Q^\gamma_{N})_{N}$ is tight, following chap. 11 of \cite{K.L.}. For that we need to show that
\begin{eqnarray*}
(1)\quad
\lim_{A\rightarrow{+\infty}}\limsup_{N\rightarrow{+\infty}}\mathbb{P}_{\nu_{\alpha}}^{\gamma}\Big(\sup_{0\leq{t}\leq{T}}||Y_{t}||_{-k}^{2}\Big)<{\infty}
\\
(2)\forall{\epsilon>{0}},\quad
\lim_{\delta\rightarrow{0}}\limsup_{N\rightarrow{+\infty}}\mathbb{P}_{\nu_{\alpha}}^{\gamma}\Big[\omega_{\delta}(Y)\geq{\epsilon}\Big]=0,
\end{eqnarray*}
where
\begin{equation*}
\omega_{\delta}(Y)=\sup_{\substack{|s-t|<\delta\\0\leq{s,t}\leq{T}}}\|Y_{t}-Y_{s}\|_{-k}.
\end{equation*}

We start by showing condition (1). For each integer $z\geq{0}$, recall that $h_{z}$ denotes the Hermite function defined at the beginning of the
second section. Denote by $M_{t}^{N,z}$ the martingale $M_{t}^{N,h_{z}}$ as defined in expression (\ref{martingaleM}).
\begin{lemma} \label{th:lemmatight}
There exists a finite constant $C(\alpha,T)$ such that for every $z\geq{0}$,
\begin{equation*}
\limsup
_{N\rightarrow{+\infty}}\mathbb{E}_{\nu_{\alpha}}^{\gamma}\Big(\sup_{0\leq{t}\leq{T}}|<Y_{t},h_{z}>|^{2}\Big)\leq{C(\alpha,T)\{<h_{z},h_{z}>\}}.
\end{equation*}
In this expression $<Y_{t},h_{z}>$ denotes the inner product of  $Y_{t}\in{\mathcal{H}_{-k}}$ and $h_{z}\in{\mathcal{H}_{k}}$.
\end{lemma}
\begin{proof}
By definition, we have that
\begin{equation*}
<Y_{t}^{N},h_{z}>=M_{t}^{N,z}+<Y_{0}^{N},h_{z}>+\int_{0}^{t}\Gamma_{1}^{h_{z}}(s)ds,
\end{equation*}
where $\Gamma_{1}^{h_{z}}(s)$ was defined in (\ref{intpartofmartingale}). To prove the Lemma, we estimate separately the
$L^{2}(\mathbb{P}_{\nu_{\alpha}}^{\gamma})$-norm of the terms on the right hand side of last equality. A simple computation shows that
\begin{equation*}
\lim_{N\rightarrow{+\infty}}\mathbb{E}_{\nu_{\alpha}}^{\gamma}(|<Y_{0}^{N},h_{z}>|)^{2}=\chi(\alpha)<h_{z},h_{z}>.
\end{equation*}
The $L^{2}(\mathbb{P}_{\nu_{\alpha}}^{\gamma})$-norm for the martingale term vanishes, combining Doob inequality with the fact that
$\mathbb{E}_{\nu_{\alpha}}^{\gamma}[(M_{t}^{N,z})^{2}]$ vanishes as $N\rightarrow{+\infty}$, for every $t\geq{0}$, namely:
\begin{equation*}
 \mathbb{E}_{\nu_{\alpha}}^{\gamma}\Big(\sup_{0\leq{t}\leq{T}}|M^{N,z}_{t}|^{2}\Big)\leq{4
\mathbb{E}_{\nu_{\alpha}}^\gamma\Big(|M^{N,z}_{T}|^{2}\Big)}.
\end{equation*}
To end, it remains to bound:
\begin{equation*}
\mathbb{E}_{\nu_{\alpha}}^{\gamma}\Big[\sup_{0\leq{t}\leq{T}}\Big(\int_{0}^{t}\Gamma_{1}^{h_{z}}(s)ds\Big)^{2}\Big].
\end{equation*}
 The idea to estimate last integral is the same as we used when analyzing the integral part of
the martingale, see the sixth section. By doing so, we have to bound
\begin{equation*}
\mathbb{E}_{\nu_{\alpha}}^{\gamma}\Big[\sup_{0\leq{t}\leq{T}}\Big(\int_{0}^{t}
\frac{N^{\gamma}}{N^{3/2}}\sum_{x\in{\mathbb{Z}}}\Delta_{N}U_{s}^{N}{h_{z}}\Big(\frac{x}{N}\Big)\bar{\eta}_{s}(x)ds\Big)^{2}\Big],
\end{equation*}
\begin{equation*}
\mathbb{E}_{\nu_{\alpha}}^{\gamma}\Big[\sup_{0\leq{t}\leq{T}}\Big(\int^{t}_{0}\frac{N^{\gamma}}{\sqrt{N}}\sum_{x\in\mathbb{Z}}\Big\{\nabla^{N}
U_{s}h_{z}\Big(\frac{x}{N}\Big)-\partial_{u}U_{s}h_{z}\Big(\frac{x}{N}\Big)\Big\}v[\eta_{s}(x)-\alpha]ds\Big)^{2}\Big],
\end{equation*}
and
\begin{equation} \label{eq:tightness1}
\mathbb{E}_{\nu_{\alpha}}^{\gamma}\Big[\sup_{0\leq{t}\leq{T}}\Big(\int_{0}^{t}
\frac{N^{\gamma}}{\sqrt{N}}\sum_{x\in{\mathbb{Z}}}\nabla^{N}{h_{z}}\Big(\frac{x}{N}\Big)\bar{\eta}_{s}(x)\bar{\eta}_{s}(x+1)ds\Big)^{2}\Big],
\end{equation}
where
\begin{equation*}
\Delta_{N}H\Big(\frac{x}{N}\Big)=N^2\Big(H\Big(\frac{x+1}{N}\Big)+ H\Big(\frac{x-1}{N}\Big)-2H\Big(\frac{x}{N}\Big)\Big).
\end{equation*}

By Schwarz inequality and since $\nu_{\alpha}$ is an invariant product measure, the first integral is bounded by
$CT^{2}N^{2\gamma-2}\frac{1}{4N}\sum_{x\in{\mathbb{Z}}}\Big(\Delta_{N}{U_{s}^{N}h_{z}}\Big(\frac{x}{N}\Big)\Big)^{2}\alpha(1-\alpha)$, which
vanishes as $N\rightarrow{+\infty}$. By Taylor expansion the second expectation vanishes.

In order to bound the last integral, we use the same idea as in the Boltzmann-Gibbs Principle, see last section. Then, we can bound the
expectation (\ref{eq:tightness1}) by
\begin{equation*}
\mathbb{E}_{\nu_{\alpha}}^{\gamma}\Big[\sup_{0\leq{t}\leq{T}}\Big(\int_{0}^{t}\frac{N^{\gamma}}{\sqrt{N}}
\sum_{j\geq{1}}\sum_{x\in{I}_{j}}\nabla^{N}h_{z}\Big(\frac{x}{N}\Big)\bar{\eta}_{s}(x)\bar{\eta}_{s}(x+1)ds\Big)^{2}\Big],
\end{equation*}
where the $I_{j}$'s are taken as in the proof of Theorem \ref{bg}, with $j$ odd for instance.

By summing and subtracting $h_{z}\Big(\frac{y_{j}}{N}\Big)$, where $y_{j}$ is one point of the interval $I_{j}$, we bound last expectation by
\begin{equation*}
2\mathbb{E}_{\nu_{\alpha}}^{\gamma}\Big[\sup_{0\leq{t}\leq{T}}\Big(\int_{0}^{t} \frac{N^{\gamma}}{\sqrt{N}}\sum_{j\geq{1}}\sum_{x\in{I_{j}}}
\Big(\nabla^{N}{h_{z}}\Big(\frac{x}{N}\Big)-\nabla^{N}{h_{z}}\Big(\frac{y_{j}}{N}\Big)\Big)\bar{\eta}_{s}(x)\bar{\eta}_{s}(x+1)ds\Big)^{2}\Big]
\end{equation*}
\begin{equation} \label{eq:tightness2}
+2\mathbb{E}_{\nu_{\alpha}}^{\gamma}\Big[\sup_{0\leq{t}\leq{T}}\Big(\int_{0}^{t}
\frac{N^{\gamma}}{\sqrt{N}}\sum_{j\geq{1}}\nabla^{N}{h_{z}}\Big(\frac{y_{j}}{N}\Big)\sum_{x\in{I_{j}}}\bar{\eta}_{s}(x)\bar{\eta}_{s}(x+1)ds\Big)^{2}\Big].
\end{equation}

By Schwarz inequality and since $\nu_{\alpha}$ is an invariant product measure, the first integral vanishes if $K$ is such that $K
N^{\gamma-1}\rightarrow{0}$ as $N\rightarrow{+\infty}$ (see Lemma \ref{lemma:lenghK1}).

To bound (\ref{eq:tightness2}), we sum and subtract inside the sum over $j$, the expectation of
$\sum_{x\in{I_{j}}}\bar{\eta}_{s}(x)\bar{\eta}_{s}(x+1)$ conditioned on the hyperplanes $M_{j}=\sigma(\sum_{x\in{I_{j}^{*}}}\eta(x))$, where
$I_{j}^{*}$ denotes the extended interval $I_{j}$. Then, we need to bound
\begin{equation*}
\mathbb{E}_{\nu_{\alpha}}^{\gamma}\Big[\sup_{0\leq{t}\leq{T}}\Big(\int_{0}^{t}\frac{N^{\gamma}}{\sqrt{N}}\sum_{j\geq{1}}
\nabla^{N}{h_{z}}\Big(\frac{y_{j}}{N}\Big)V_{j}(\eta_{s})ds\Big)^{2}\Big]
\end{equation*}
and
\begin{equation*}
\mathbb{E}_{\nu_{\alpha}}^{\gamma}\Big[\sup_{0\leq{t}\leq{T}}\Big(\int_{0}^{t}
\frac{N^{\gamma}}{\sqrt{N}}\sum_{j\geq{1}}\nabla^{N}{h_{z}}\Big(\frac{y_{j}}{N}\Big)E\Big(\sum_{x\in{I_{j}}}\bar{\eta}_{s}(x)\bar{\eta}_{s}(x+1)|M_{j}\Big)
ds\Big)^{2}\Big].
\end{equation*}
where
\begin{equation*}
V_{j}(\eta)=\sum_{x\in{I_{j}}}\bar{\eta}_{s}(x)\bar{\eta}_{s}(x+1)- E\Big(\sum_{x\in{I_{j}}}\bar{\eta}_{s}(x)\bar{\eta}_{s}(x+1)|M_{j}\Big).
\end{equation*}
By Lemma 4.3 of \cite{C.L.O}, the first integral is bounded by
\begin{equation*}
C_{0}\int_{0}^{T}\Big|\Big|\frac{N^{\gamma}}{\sqrt{N}}\sum_{j\geq{1}} \nabla^{N}{h_{z}}\Big(\frac{y_{j}}{N}\Big)V_{j}\Big|\Big|_{-1}^{2}ds,
\end{equation*}
where $C_{0}$ is a constant. To bound this $H_{-1}$-norm one can follow the same computations as in Lemma (\ref{lm:lengthK}), and it is not hard
to show that it vanishes for $K$, such that $K^{2}N^{\gamma-1}\rightarrow{0}$ as $N\rightarrow{+\infty}$. To bound the other integral, we
proceed in the same lines as in the Boltzmann-Gibbs Principle, which concludes the proof.
\end{proof}

\begin{corollary}
For each $k>1$
\begin{equation*}
(a)\quad \limsup_{N\rightarrow{+\infty}}\mathbb{E}_{\nu_{\alpha}}^{\gamma}\Big(\sup_{0\leq{t}\leq{T}}||Y_{t}||_{-k}^{2}\Big)<{\infty}
\end{equation*}
\begin{equation*}
(b)\quad
\lim_{n\rightarrow{+\infty}}\limsup_{N\rightarrow{+\infty}}\mathbb{E}_{\nu_{\alpha}}^{\gamma}\Big[\sup_{0\leq{t}\leq{T}}\sum_{|z|\geq{n}}(<Y_{t},h_{z}>)^{2}\gamma_{z}^{-k}\Big]=0.
\end{equation*}
\end{corollary}
\begin{proof}
Recall the definition of $\mathcal{H}_{k}$ and the inner product $<,>_{k}$ at the beginning of the second section. Since
$<f,g>_{k}=\sum_{z\in{\mathbb{Z}}}<f,h_{z}><g,h_{z}>\gamma_{z}^{-k}$, then
\begin{equation*}
\limsup_{N\rightarrow{+\infty}}\mathbb{E}_{\nu_{\alpha}}^{\gamma}\Big(\sup_{0\leq{t}\leq{T}}||Y_{t}||_{-k}^{2}\Big)\leq{\limsup_{N\rightarrow{+\infty}}
\sum_{z\in{\mathbb{Z}}}\gamma_{z}^{-k}\mathbb{E}_{\nu_{\alpha}}^{\gamma}\Big(\sup_{0\leq{t}\leq{T}}<Y_{t},h_{z}>^{2}\Big)}.
\end{equation*}
 and by the previous Lemma it is bounded by $C(\alpha,T)\sum_{z\in{\mathbb{Z}}}\gamma_{z}^{-k}$,
which is finite as long as $k>1$. The assertion (b) follows by the same argument.
\end{proof}
We note that this is the place where we need the restriction $k>1$ in order to have the density fluctuation field well defined in
$\mathcal{H}_{-k}$.

By Chebychev's inequality the condition (1) is a consequence of (a) in the previous Corollary. So, in order to prove that
$(Q_{N}^{\gamma})_{N\geq{1}}$ is tight we only have to show statement (2). In view of (b), this follows from the next Lemma.
\begin{lemma}
For every $n\in{\mathbb{N}}$ and every $\epsilon>{0}$,
\begin{equation*}
\lim_{\delta\rightarrow{0}}\limsup_{N\rightarrow{+\infty}}\mathbb{P}_{\nu_{\alpha}}^{\gamma}\Big[\sup_{\substack{|s-t|<\delta\\0\leq{s,t}\leq{T}}}\quad
\sum_{|z|\leq{n}}(<Y_{t}-Y_{s},h_{z}>)^{2}\gamma_{z}^{-k}>\epsilon\Big]=0.
\end{equation*}
\end{lemma}
To prove this Lemma it is enough to show that
\begin{equation*}
\lim_{\delta\rightarrow{0}}\limsup_{N\rightarrow{+\infty}}\mathbb{P}_{\nu_{\alpha}}^{\gamma}\Big[\sup_{\substack{|s-t|<\delta\\0\leq{s,t}\leq{T}}}\quad
(<Y_{t}-Y_{s},h_{z}>)^{2}>\epsilon\Big]=0.
\end{equation*}
for every $z\in{\mathbb{Z}}$ and $\epsilon>0$, which is a consequence of the next two results.

\begin{lemma}
Fix a function $H\in{S(\mathbb{R})}$. For every $\epsilon>0$
\begin{equation*}
\lim_{\delta\rightarrow{0}}\limsup_{N\rightarrow{+\infty}}\mathbb{P}_{\nu_{\alpha}}^{\gamma}\Big[\sup_{\substack{|s-t|<\delta\\0\leq{s,t}\leq{T}}}\quad
|M_{t}^{N,H}-M_{s}^{N,H}|>\epsilon\Big]=0.
\end{equation*}
\end{lemma}
\begin{proof}
Denote by $\omega'_{\delta}(M^{H})$ the modified modulus of continuity defined by
\begin{equation*}
\omega'_{\delta}(M^{N,H})=\inf_{\substack{\{t_{i}\}}}\quad\max_{\substack{0\leq{i}\leq{r}}}\quad\sup_{\substack{t_{i}\leq{s}<{t}\leq{t_{i+1}}}}|M_{t}^{N,H}
-M_{s}^{N,H}|,
\end{equation*}
where the infimum is taken over all partitions of $[0,T]$ such that $0=t_{0}<t_{1}<...<t_{r}=T$ where $t_{i+1}-t_{i}>\delta$ for
$0\leq{i}\leq{r}$.

Since
\begin{equation*}
\omega_{\delta}(M^{N,H})\leq{2\omega'_{\delta}(M^{N,H})+\sup_{\substack{t}}|M_{t}^{N,H}-M_{t_{-}}^{N,H}|}
\end{equation*}
and
\begin{equation*}
\sup_{\substack{t}}\Big|M_{t}^{N,H}-M_{t_{-}}^{N,H}\Big|=\sup_{\substack{t}}\Big|
<Y_{t}^{N},H>-<Y_{t_{-}}^{N},H>\Big|\leq{\frac{||\nabla{H}||_{\infty}}{N^{1+\frac{1}{2}}}},
\end{equation*}
the proof ends if one shows that
\begin{equation*}
\lim_{\delta\rightarrow{0}}\limsup_{N\rightarrow{+\infty}}\mathbb{P}_{\nu_{\alpha}}^{\gamma}\Big[\omega'_{\delta}(M^{N,H})>\epsilon\Big]=0
\end{equation*}
for every $\epsilon>0$. By the Aldous criterium, see for example Proposition 4.1.6 of \cite{K.L.} it is enough to show that:
\begin{equation*}
\lim_{\delta\rightarrow{0}}\limsup_{N\rightarrow{+\infty}}\sup_{\substack{\tau\in{\mathfrak
{T}_{\tau}}\\0\leq{\theta}\leq{\delta}}}\mathbb{P}_{\nu_{\alpha}}^{\gamma}\Big[|M_{\tau+\theta}^{N,H}-M_{\tau}^{N,H}|>\epsilon\Big]=0
\end{equation*}
for every $\epsilon>0$. Here $\mathfrak {T}_{\tau}$ denotes the family of all stopping times bounded by $T$ with respect to the canonical
filtration. By Chebychev´s inequality, the Optional Sampling Theorem and expression (\ref{quadraticlonger})
 the result follows.
\end{proof}
\begin{lemma}
Fix $H\in{S(\mathbb{R})}$. For every $\epsilon>0$
\begin{equation*}
\lim_{\delta\rightarrow{0}}\limsup_{N\rightarrow{+\infty}}
\mathbb{P}_{\nu_{\alpha}}^{\gamma}\Big[\sup_{\substack{|s-t|<\delta\\0\leq{s,t}\leq{T}}}\Big|
\int_{s}^{t}\Gamma_{1}^{H}(r)dr\Big|>\epsilon\Big]=0
\end{equation*}
\end{lemma}
\begin{proof}
By using the explicit knowledge of $\Gamma_{1}^{H}(r)$, the decomposition of the instantaneous current (\ref{eq:decompcurrent}) and similar
computations as the ones performed when analyzing the integral part of the martingale $M_{t}^{N,H}$, we just need to bound:
\begin{equation*}
\mathbb{P}_{\nu_{\alpha}}^{\gamma}\Big[\sup_{\substack{|s-t|<\delta\\0\leq{s,t}\leq{T}}}\Big|
\int_{s}^{t}\frac{N^{\gamma}}{\sqrt{N}}\sum_{x\in{\mathbb{Z}}}\nabla^{N}U_{r}^{N}H\Big(\frac{x}{N}\Big)
\bar{\eta}_{r}(x)\bar{\eta}_{r}(x+1)dr\Big|>\epsilon\Big].
\end{equation*}
Dividing the interval $[0,T]$ in small intervals of length $\delta$, last probability is bounded by
\begin{equation*}
\frac{T}{\delta}\mathbb{P}_{\nu_{\alpha}}^{\gamma}\Big[\sup_{0\leq{t\leq{\delta}}}\Big|
\int_{0}^{t}\frac{N^{\gamma}}{\sqrt{N}}\sum_{x\in{\mathbb{Z}}}\nabla^{N}U_{r}^{N}H\Big(\frac{x}{N}\Big)\bar{\eta}_{s}(x)\bar{\eta}_{s}(x+1)dr\Big|>\frac{\epsilon}{2}\Big].
\end{equation*}
Using Chebychev's inequality, last probability is bounded by the expectation that appeared at the end of the proof of Lemma
(\ref{th:lemmatight}) which we showed to vanish as $N\rightarrow{+\infty}$.
\end{proof}

\section{Dependence on the initial configuration for the longer time scale}

We start by considering the case $\alpha=1/2$ which implies that $v=0$. In this case, we can define ( as in the hyperbolic scaling) for a site
$x$, the current over the fixed bond $[x,x+1]$ denoted by $J^{N,\gamma}_{x,x+1}(t)$, as the total number of jumps from the site $x$ to the site
$x+1$ minus the total number of jumps from the site $x+1$ to the site $x$ during the time interval $[0,tN^{1+\gamma}]$.

In this particular case, the density fluctuation field at time $t$ is the same as at time $0$. As a consequence, the current through $[x,x+1]$
converges to $0$ in the $L^{2}(\mathbb{P}_{\nu_{\alpha}}^{\gamma})$-norm:

\begin{proposition} \label{current dependence higher}
Fix $t\geq{0}$, a site $x\in{\mathbb{Z}}$ and $\gamma<1/3$. Then,
\begin{equation*}
\lim_{N\rightarrow{+\infty}}\mathbb{E}_{\nu_{\alpha}}^{\gamma}\Big[\frac{\bar{J}^{N,\gamma}_{x,x+1}(t)}{\sqrt{N}}\Big]^{2}=0.
\end{equation*}
\end{proposition}

The idea of the proof is the same as the one used in the hyperbolic scaling, and it relies on the following result:
\begin{proposition} \label{prop:2}
For every $t\geq{0}$ and $\gamma<1/3$:
\begin{equation*}
\lim_{n\rightarrow{+\infty}}\mathbb{E}_{\nu_{\alpha}}^{\gamma}\Big[\frac{\bar{J}^{N,\gamma}_{-1,0}(t)}{\sqrt{N}}-(Y_{t}^{N,\gamma}(G_{n})-Y_{0}^{N,\gamma}(G_{n}))\Big]^{2}=0,
\end{equation*}
uniformly over $N$.
\end{proposition}
\begin{proof}
Recall the proof of Proposition \ref{prop:1}. There is only a slight difference that we need to remark. In this case, the expression
(\ref{eq:current}) in the proof of that Proposition becomes equal to
\begin{equation*}
\frac{1}{\sqrt{N}}\sum_{x=1}^{Nn}\frac{1}{Nn}M^{N,\gamma}_{x-1,x}(t)
+\frac{N^{\gamma}}{\sqrt{N}}\int_{0}^{t}\frac{1}{n}\sum_{x=1}^{Nn}[W_{x,x+1}(\eta_{s})-(p-q)\chi(\alpha)]ds,
\end{equation*}
where $M^{N,\gamma}_{x-1,x}(t)$ denotes the martingale associated to the current through the bond $[x-1,x]$.

 Estimating the quadratic variation of the martingale $M^{N,\gamma}_{x-1,x}(t)$ by $N^{1+\gamma}t$, the $L^{2}(\mathbb{P}_{\nu_{\alpha}}^{\gamma})$-norm of the
 martingale term in last equality, is bounded by $\frac{N^{1+\gamma}t}{n}$ which vanishes as $n\rightarrow{+\infty}$.
To bound the integral term, using the decomposition of the instantaneous current (\ref{eq:decompcurrent}), it is enough to bound
\begin{equation*}
\mathbb{E}_{\nu_{\alpha}}^{\gamma}\Big[\frac{N^{\gamma}}{\sqrt{N}}\int_{0}^{t}\frac{1}{n}\sum_{x=0}^{Nn-1}\bar{\eta}_{s}(x)\bar{\eta}_{s}(x+1)ds\Big]^{2}.
\end{equation*}
Using the inequality $(x+y)^2\leq{2x^2+2y^2}$, last expectation is bounded by
\begin{equation*}
2\mathbb{E}_{\nu_{\alpha}}^{\gamma}\Big[\frac{(Nn)^{\gamma}}{\sqrt{nN}}\int_{0}^{t}\sum_{x=0}^{Nn}\bar{\eta}_{s}(x)\bar{\eta}_{s}(x+1)ds\Big]^{2}+
2\mathbb{E}_{\nu_{\alpha}}^{\gamma}\Big[\frac{N^{\gamma}}{\sqrt{N}}\int_{0}^{t}\frac{1}{n}\bar{\eta}_{s}(Nn-1)\bar{\eta}_{s}(Nn)ds\Big]^{2}.
\end{equation*}

Recall the proof of the Boltzmann-Gibbs Principle when applied to the function $H(u)=1_{[0,1]}(u)$, which gives us that the expectation on the
left hand side of last expression vanishes as $n\rightarrow{+\infty}$ uniformly over $N$. By Schwarz inequality and since $\nu_{\alpha}$ is an
invariant product measure, the other term vanishes as $n\rightarrow{+\infty}$, which concludes the proof.
\end{proof}
Last result is stated for the bond $[-1,0]$ but for $[x,x+1]$ a similar statement holds.

Consider now the case $\alpha\neq{1/2}$. In this case, by the definition of the density fluctuation field (see (\ref{eq:densfieldlongscale})),
as time is going by the position of the particles start to change. So, if there is initially a particle at site $x$ and if it does not move,
then at time $t$, its position is the site $x+[vtN^{1+\gamma}]$, that we denote by $y^{x}_{t}$. By this reason, we cannot consider any longer
the current through a fixed bound, but we must consider the current through a bond that depends on time.

Let $J^{N,\gamma}_{y^{x}_{t}}(t)$ be the current trough the bond $[y^{x}_{t},y^{x}_{t}+1]$, defined as the number of particles that jump from
$y^{x}_{t}$ to $y^{x}_{t}+1$, minus the number of particles that jump from $y^{x}_{t}+1$ to $y^{x}_{t}$, from time $0$ to $tN^{1+\gamma}$.
Formally we have that:
\begin{equation*}
J^{N,\gamma}_{y^{x}_{t}}(t)=\sum_{y\geq{1}}\Big(\eta_{t}(y+y^{x}_{t})-\eta_{0}(y+x)\Big).
\end{equation*}
As a consequence, it holds that:
\begin{proposition} \label{current dependence higher}
Fix $t\geq{0}$, a site $x\in{\mathbb{Z}}$ and $\gamma<1/3$. Then,
\begin{equation*}
\lim_{N\rightarrow{+\infty}}\mathbb{E}_{\nu_{\alpha}}^{\gamma}\Big[\frac{\bar{J}^{N,\gamma}_{y^{x}_{t}}(t)}{\sqrt{N}}\Big]^{2}=0.
\end{equation*}
\end{proposition}
As in the hyperbolic scaling, this last results is a consequence of the following:
\begin{proposition}
For every $t\geq{0}$ and $\gamma<1/3$:
\begin{equation*}
\lim_{n\rightarrow{+\infty}}\mathbb{E}_{\nu_{\alpha}}^{\gamma}\Big[\frac{\bar{J}^{N,\gamma}_{y^{x}_{t}}(t)}{\sqrt{N}}-(Y_{t}^{N,\gamma}(G_{n})-Y_{0}^{N,\gamma}(G_{n}))\Big]^{2}=0,
\end{equation*}
uniformly over $N$.
\end{proposition}
\begin{proof}
Recall the proof of Proposition \ref{prop:2}. The martingale associated to $J^{N,\gamma}_{y^{x}_{t}}(t)$ is now given by
\begin{equation*}
M^{N,\gamma}_{x}(t)=J^{N,\gamma}_{y^{x}_{t}}(t)-\int_{0}^{t}
\Big\{N^{1+\gamma}W_{y^{x}_{s}}(\eta_{s})+\partial_{s}J^{N,\gamma}_{y^{x}_{s}}(s)\Big \}ds,
\end{equation*}
where $W_{y^{x}_{s}}(\eta)$ denotes the instantaneous current through the bond $[y^{x}_{t},y^{x}_{t}+1]$. Since
$\partial_{s}J^{N,\gamma}_{y^{x}_{t}}(s)=-vN^{1+\gamma}\eta_{s}(y^{x}_{s})$ and repeating the same arguments as in the proof of Proposition
\ref{prop:2}, the result follows.
\end{proof}

\quad\

\demons\\ In this case there is a relation between the position of the Tagged particle and the current through the bond
$[y^{-1}_{t},y^{-1}_{t}+1]$ and the density of particles, which is given by:
\begin{equation*}
\Big\{X_{tN^{1+\gamma}}\geq{a}\Big\}=\Big\{J^{N,\gamma}_{y^{-1}_{t}}(t)\geq{\sum_{x=vtN^{1+\gamma}}^{a-1}\eta_{t}(x)}\Big\}.
\end{equation*}
Repeating the same computations as in the proof of Corollary \ref{tagged dependence}, using the fact that
$\mathbb{E}_{\nu_{\alpha}}^{\gamma}[J^{N,\gamma}_{y^{-1}_{t}}(t)]=(p-q)\alpha^2 tN^{1+\gamma}$; that
$\frac{\bar{J}^{N,\gamma}_{y^{-1}_{t}}(t)}{\sqrt{N}}$ converges to $0$ in the $L^{2}(\mathbb{P}_{\nu_{\alpha}}^{\gamma})$-norm; that
\begin{equation*}
\lim_{N\rightarrow{+\infty}}\mathbb{E}_{\nu_{\alpha}}^{\gamma}\Big[Y_{t}^{N,\gamma}(H)-Y_{0}^{N,\gamma}(H)\Big]^{2}=0
\end{equation*}
for every $H\in{S(\mathbb{R})}$ and also that
\begin{equation*}
\mathbb{E}_{\nu_{\alpha}}^{\gamma}\Big[\frac{1}{\sqrt{N}}\sum_{x=1+v_{t}N^{1+\gamma}}^
{a\sqrt{N}-1+v_{t}N^{1+\gamma}+Z^{N,\gamma}_t}\bar{\eta}_{t}(x)\Big]^2=O(N^{-1/2}),
\end{equation*}
where $Z^{N,\gamma}_t$
\begin{equation*}
Z^{N,\gamma}_t=-\sum_{x=1+vtN^{1+\gamma}}^{v_{t}N^{1+\gamma}}\bar{\eta}_{t}(x)/\alpha
\end{equation*}
 the result follows. \cqd
\\
\textbf{Ackowlodgements:} \textit{I am deeply indebted to my Phd advisor, Claudio Landim, for suggesting this problem and for valuable
discussions. Also, I would like to express my gratitude to Sunder Sethuraman for his help on explaining some topics related to the Theorem (2.1)
of \cite{S.}. Finally, I want to thank Johel Beltran for stimulating discussions.}

\end{document}